\documentclass[global]{svjour}

\usepackage{amsfonts}
\usepackage{amsmath}
\usepackage{amssymb}
\usepackage{fullpage}
\usepackage{graphicx,rotating}

\newcommand{\mbb}[1]{\mathbb{#1}}
\newcommand{\mb}[1]{\mathbf{#1}}
\newcommand{\wt}[1]{\widetilde{#1}}

\journalname{Electronic Journal of Probability}

\begin{document}

\title{Record indices and age-ordered frequencies in Exchangeable Gibbs
Partitions}
\author{\large {Robert C. Griffiths \and Dario Span\`{o}}\thanks{Research supported in part by EPSRC grant GR/T21783/01.}}
\institute {\large University of Oxford
\\ Department of Statistics, 1 South Parks Road, Oxford OX1 3TG, UK\\
\email{\large \{griff,spano\}@stats.ox.ac.uk}}
\date{\large \today}

\mail{Dario Span\`{o}, Department of Statistics, 1 South Parks Road,
Oxford OX1 3TG, United Kingdom.}

\maketitle

\abstract{The frequencies $X_{1},X_{2},\ldots$ of an exchangeable
Gibbs random partition $\Pi$ of $\mbb{N}=\{1,2,\ldots\}$ (Gnedin
and Pitman (2006)) are considered in their \emph{age-order}, i.e.
their size-biased order. We study their dependence on the sequence
$i_1,i_2,\ldots$ of least elements of the blocks of $\Pi$. In
particular, conditioning on  $1=i_{1}<i_{2}<\ldots$, a
representation is shown to be
 $$X_{j}=\xi_{j-1}\prod_{i=j}^{\infty}(1-\xi_{i}) \hspace{2cm}j=1,2,\ldots$$
 where $\{\xi_{j}:j=1,2,\ldots\}$ is a sequence of independent Beta random variables.
Sequences with such a product form are called \emph{neutral to the
left}. We show that the property of conditional left-neutrality in
fact characterizes the Gibbs family among all exchangeable
partitions,
   and leads to further interesting results on: (i) the conditional Mellin transform of
 $X_k$, given $i_{k}$, and (ii) the conditional distribution of the first $k$ normalized
 frequencies, given $\sum_{j=1}^{k}X_j$ and $i_k$; the latter turns out to be a mixture of Dirichlet distributions.
Many of the mentioned representations are extensions of Griffiths
and Lessard (2005) results on Ewens' partitions.


\newpage

\section{Introduction.}

A random partition of $[n]=\{1,\ldots,n\}$ is a random collection
$\Pi_{n}=\{\Pi_{n1},\ldots,\Pi_{nk}\}$ of disjoint nonempty
subsets of $[n]$ whose union is $[n]$. The classes of $\Pi_{n}$
are conventionally ordered by their least elements
$1=i_{1}<i_{2}<\ldots<i_{k}\leq n$. We call $\{i_{j}\}$  the
sequence of \emph{record indices} of $\Pi_{n}$, and define the
\emph{age-ordered
 frequencies} of $\Pi_{n}$ to be the vector
 $\mb{n}=(n_{1},\ldots,n_{k})$ such that $n_{j}$ is the cardinality
 of $\Pi_{nj}$.
  Consistent
Markov partitions $\Pi=(\Pi_{1},\Pi_{2},\ldots)$ can be generated by
a set of predictive distributions specifying, for each $n$, how
$\Pi_{n+1}$ is likely to extend $\Pi_{n}$, that is: given $\Pi_{n}$,
a conditional probability is assigned for the integer $(n+1)$ to
 join any particular class of $\Pi_n$ or to start a new class.

\noindent We consider a family of consistent random partitions
studied by Gnedin and Pitman \cite{GP} which can be defined by the
following prediction rule: (i) set $\Pi_{1}=(\{1\})$; (ii) for each
$n\geq 1$, conditional on $\Pi_{n}=(\pi_{n1},\ldots,\pi_{nk})$, the
probability that $(n+1)$ starts a new class is

\begin{equation}
\frac{V_{n+1,k+1}}{V_{n,k}}\label{gn}
\end{equation}
otherwise, if $n_{j}$ is the cardinality of $\pi_{nj}$
$(j=1,\ldots,k)$, the probability that $(n+1)$ falls in the $j$-th
\lq\lq old\rq\rq class $\pi_{nj}$ is
\begin{equation}
\frac{n_{j}-\alpha}{n-\alpha
k}\left(1-\frac{V_{n+1,k+1}}{V_{n,k}}\right),\label{go}
\end{equation} for some $\alpha\in(-\infty,1]$ and a sequence of coefficients
$V=(V_{n,k}:k\leq n=1,2,\ldots)$ satisfying the recursion:
 \begin{equation}
V_{1,1}=1;\hspace{2cm}V_{n,k}=(n-\alpha k)V_{n+1,k}+V_{n+1,k+1}.
\label{v}
\end{equation}
Every partition $\Pi$ of $\mathbb{N}$ so generated is called an
\emph{exchangeable Gibbs partition with parameters $(\alpha,V)$}
(EGP$(\alpha,V)$), where exchangeable means that, for every $n$, the
distribution of $\Pi_{n}$ is a symmetric function of the vector
$\mb{n}=(n_{1},\ldots,n_{k})$ of its frequencies (\cite{P95}) (see
section \ref{o} below). Actually, the whole family of EGPs, treated
in \cite{GP} includes also the value $\alpha=-\infty$, for which the
definition (\ref{gn})-(\ref{go}) should be modified; this case will
not be treated in the present paper.

\noindent A special subfamily of EGPs is \emph{Pitman's
two-parameter family}, for which $V$ is given by
\begin{equation}
V^{(\alpha,\theta)}_{n,k}=\frac{\prod_{j=1}^{k}(\theta+\alpha(j-1))}{\theta_{(n)}}\label{p2p}
\end{equation}
where either $\alpha\in [0,1]$ and $\theta\geq -\alpha$ or
$\alpha<0$ and $\theta=m|\alpha|$ for some integer $m$. Here and in
the following sections, $a_{(x)}$ will denote the generalized
increasing factorial i.e. $a_{(x)}=\Gamma(a+x)/\Gamma(a)$, where
$\Gamma(\cdot)$ is the Gamma function.

\noindent Pitman's family is characterized as the unique class of
EGPs with $V$-coefficients of the form
$$V_{n,k}=\frac{V^{*}_{k}}{c_{n}}$$ for some sequence of constants
$(c_{n})$ (\cite{GP}, Corollary 4). If we let $\alpha=0$ in
(\ref{p2p}), we obtain the well known \emph{Ewens' partition} for
which:
\begin{equation}
V^{(0,\theta)}_{n,k}=\frac{\theta^{k}}{\theta_{(n)}}\label{ew}.
\end{equation}
Ewens' family arose in the context of Population Genetics to
describe the properties of a population of genes under the
so-called infinitely-many-alleles model with parent-independent
mutation (see e.g. \cite{W84}, \cite{K82}) and became a paradigm
for the modern developments of a theory of exchangeable random
partitions (\cite{A85}, \cite{P95}, \cite{G97}).

\ \\
\noindent For every fixed $\alpha$, the set of all EGP$(\alpha,V)$
forms a convex set; Gnedin and Pitman proved it and gave a
complete description of the extreme points (\cite{GP}, Theorem
12). It turns out, in particular, that for every $\alpha\leq 0$,
the extreme set is given by Pitman's two-parameter family. For
each $\alpha\in(0,1)$, the extreme points are all partitions of
the so-called Poisson-Kingman type with parameters $(\alpha,s), \
s>0$, whose $V$-coefficients are given by:

\begin{equation}
V_{n,k}(s)=\alpha^{k}s^{n/\alpha}G_{\alpha}(n-\alpha
k,s^{-1/\alpha}), \label{pkv}
\end{equation}
with
\begin{equation}G_{\alpha}(q,t):=\frac{1}{\Gamma(q)f_{\alpha}(t)}\int_{0}^{t}f_{\alpha}(t-v)v^{q-1}dv,\label{galpha}\end{equation}
and $f_{\alpha}$ is an $\alpha$-stable density (\cite{P02},
Theorem 4.5). The partition induced by (\ref{pkv}) has
 limit
frequencies (ranked in a decreasing order) equal in distribution to
the jump-sizes of the process $(S_{t}/S_{1}:t\in[0,1])$,
{conditioned on $S_{1}=s$}, where $(S_{t}:t>0)$ is a stable
subordinator
 with density $f_{\alpha}$. The parameter $s$ has the interpretation
as the (a.s.) limit of the ratio $K_{n}/n^{\alpha}$ as
$n\rightarrow\infty$, where $K_{n}$ is the number of classes in the
partition $\Pi_{n}$ generated via $V_{n,k}(s)$. $K_n$ is shown in
\cite{GP} to play a central role in determining the extreme set of
$V_{n,k}$ for every $\alpha$; the distribution of $K_n$, for every
$n$, turns out to be of the form
\begin{equation}
\mbb{P}(K_n=k)=V_{n,k}{n\brack k}_{\alpha}.
\label{knpr}\end{equation}\\
where  ${n\brack k}_{\alpha}$ are generalized Stirling numbers,
defined as the coefficients of $x^{n}$ in
$$\frac{n!}{\alpha^{k}k!}(1-(1-x)^{\alpha})^{k}$$
(see \cite{GP} and reference therein). As $n\rightarrow\infty$,
$K_n$ behaves differently for different choices of the parameter
$\alpha$: almost surely it will be finite for $\alpha<0$, $K_n\sim
S\log n$ for $\alpha=0$ and $K_n\sim S n^{\alpha}$ for positive
$\alpha$, for some positive random variable $S$.
\ \\

\noindent In this paper we want to study how the distribution of
the limit age-ordered frequencies
$X_{j}=\lim_{n\rightarrow\infty}n_{j}/n$ $(j=1,2,\ldots)$ in an
Exchangeable Gibbs partition depends on its record indices
${\mb{i}}=(1=i_{1}<i_{2}<\ldots)$.
 To this purpose, we adopt a combinatorial
approach proposed by Griffiths and Lessard \cite{GL} to study the
distribution of the age-ordered allele frequencies
$X_{1},X_{2},\ldots$ in a population corresponding to the
so-called Coalescent process with mutation (see e.g. \cite{W84}),
whose equilibrium distribution is given by Ewens' partition
(\ref{ew}), for some mutation parameter $\theta>0$. In such a
context, the record index $i_j$ has the interpretation as the
number of ancestral lineages surviving back in the past, just
before the last gene of the $j$-th oldest type, observed in the
current generation, is lost by mutation.

\noindent Following Griffiths and Lessard's steps we will (i) find,
for every $n$, the distribution of the age ordered frequencies
$\mb{n}=(n_{1},\ldots,n_{k})$, conditional on the record indices
$\mb{i}_{n}=(1=i_{1}<i_{2}<\ldots<i_{k})$ of $\Pi_{n}$, as well as
the distribution of $\mb{i}_{n}$; (ii) take their limits as
$n\rightarrow\infty$; (iii) for $m=1,2,\ldots,$ describe the
distribution of the $m$-th age-ordered frequency conditional on
$i_m$ alone. We will follow such steps, respectively, in sections
\ref{*}, \ref{**}, \ref{***}. In addition, we will derive in section
\ref{****} a representation for the distribution of the first $k$
age-ordered frequencies, conditional on their cumulative sum and on
$i_k$.

\noindent In our investigation of EGPs, the key result is relative
to the step (ii), stated in Proposition \ref{cgem}, where we find
that, conditional on $\mb{i}=(1,i_{2},\ldots)$, for every
$j=1,2,\ldots$,

\begin{equation}
X_{j}|\mb{i}\overset{d}{=}\xi_{j-1}\prod_{i=j}^{\infty}(1-\xi_{i}),
\label{prepro}
\end{equation}
almost surely, for an independent sequence
$(\xi_{0},\xi_{1},\ldots)\in[0,1]^{\infty}$ such that $\xi_{0}\equiv
1$ and $\xi_{m}$ has a Beta density with parameters
$(1-\alpha,i_{m+1}-\alpha m-1)$ for each $m\geq 1$. \noindent The
representation (\ref{prepro}) does not depend on $V$. The parameter
$V$ affects only the distribution of the record indices
$\mb{i}=(i_{1},i_{2},\ldots)$ which is a non-homogenous Markov
chain, starting at $i_{1}\equiv 1$, with transition probabilities

\begin{equation}
P_{j}(i_{j+1}|i_{j})=(i_{j}-\alpha
j)_{(i_{j+1}-i_{j}-1)}\frac{V_{i_{j+1},j+1}}{V_{i_{j},j}},\hspace{1.5cm}j\geq
1. \label{itr}
\end{equation}
The representation (\ref{prepro})-(\ref{itr}) extends Griffiths and
Lessard's result on Ewens' partitions (\cite{GL}, (29)), recovered
just by letting $\alpha=0$.

\noindent In section \ref{**} we stress the connection between the
representation (\ref{prepro}) and a wide class of random discrete
distributions, known in the literature of Bayesian Nonparametric
Statistics as Neutral to the Left (NTL) processes (\cite{DK99},
\cite{BW}) and use such a connection to show that the structure
(\ref{prepro}) with independent $\{\xi_j\}$ actually characterizes
EGP's among all exchangeable partitions of $\mbb{N}$.

 \noindent The representation (\ref{prepro}) is useful to find the
 moments of
both $X_{j}$ and $\sum_{i=1}^{j}X_{i}$, conditional on the $j$-th
record index $i_{j}$ alone, as shown in section \ref{***}. In the
same section a recursive formula is found for the Mellin transform
of both random quantities, in terms of the Mellin transform of the
size-biased pick $X_1$. In section \ref{***i}, an alternative
description  of the conditional distribution of $-\log X_j$ given
$i_j$ leads to a characterization of Ewens-Kingman's partitions as
the only class for which $-\log X_j$ can be expressed as an
infinite sum of independent random variables.

\noindent Finally in section \ref{****} we obtain an expression for
the density of the first $k$ age-ordered frequencies
$X_1,\ldots,X_k$, conditional on $\sum_{i=1}^{k}X_{i}$ and $i_k$,
 as a mixture of Dirichlet distributions on the
$(k-1)$-dimensional simplex $(k=1,2,\ldots)$. Such a result leads
to a self-contained proof for the marginal distribution of $i_k$,
whose formula is closely related to Gnedin and Pitman's  result
(\ref{knpr}).\\

\noindent As a completion to our results, it should be noticed
that a representation of the \emph{unconditional} distribution of
the age-ordered frequencies of an EGP can be derived as a mixture
of the age-ordered distributions of their extreme points, which
are known: for $\alpha\leq 0$, the extreme age-ordered
distribution is the celebrated two-parameter GEM distribution
(\cite{P95},\cite{P96}), for which:
\begin{equation}
X_j=B_j\prod_{i=1}^{j-1}(1-B_j),\ \ j=1,2,\ldots \label{2gem}
\end{equation}
for a sequence $(B_j:j=1,2,\ldots)$ of independent Beta random
variables with parameters, respectively
$\{(1-\alpha,\theta+j\alpha):j=1,2,\ldots\}$. Such a representation
reflects a property of \emph{right-neutrality}, which in a sense is
the inverse of $(\ref{prepro})$, as it will be clear in section
\ref{**}. When $\alpha$ is strictly positive, the structure of the
age-ordered frequencies in the extreme points lose such a simple
structure. A description is available in \cite{PPY}.

\noindent We want to embed Griffiths and Lessard's method in the
general setting of Pitman's theory of exchangeable and partially
exchangeable random partitions, for which our main reference is
\cite{P95}. Pitman's theory will be summarized in section \ref{o}
. The key role played by record indices in the study of random
partitions has been emphasized by several authors, among which
Kerov \cite{K2}, Kerov and Tsilevich \cite{KT1}, and more recently
by Gnedin \cite{G06}, and  Nacu \cite{N06}, who showed that the
law of a partially exchangeable random partition is completely
determined by that of its record indices. We are indebted to an
anonymous referee for signalling the last two references, whose
findings have an intrinsic connection with many formulae in our
section \ref{o}.

\section{Exchangeable and partially exchangeable random partitions.}
\label{o}

\noindent We complete the introductory part with a short review of
Pitman's theory of exchangeable and partially exchangeable random
partitions, and stress the connection with the distribution of their
record indices. For more details we refer the reader to \cite{P95}
and \cite{P02} and reference therein. Let $\mu$ be a distribution on
$\Delta=\{x=(x_{1},x_{2},\ldots)\in[0,1]^{\infty}:|x|\leq 1\}$,
 endowed with a Borel sigma-field. Consider the function:

\begin{equation}
q_{\mu}(n_{1},\ldots,n_{k})=\int_{\Delta}\left(\prod_{j=1}^{k}p_{j}^{n_{j}-1}\right)\prod_{1}^{k-1}\left(1-\sum_{i=1}^{j}p_{i}\right)\mu(dp).
\label{eppf}
\end{equation}
The function $q_{\mu}$ is called the \emph{Partially exchangeable
probability function (PEPF)} of $\mu$, and has the interpretation as
the probability distribution of a random partition
$\Pi_{n}=(\Pi_{n1},\ldots,\Pi_{nk})$, for which a sufficient
statistic is given by its age-ordered frequencies
$(n_{1},\ldots,n_{k})$, that is:
$$\mbb{P}_{\mu}(\Pi_{n}=(\pi_{n1},\ldots,\pi_{nk}))=q_{\mu}(n_{1},\ldots,n_{k})$$
for every partition $(\pi_{n1},\ldots,\pi_{nk})$ such that
$|\pi_{nj}|=n_{j}\ (j=1,\ldots,k\leq n).$

\noindent If $q_{\mu}(n_{1},\ldots,n_{k})$ is symmetric with respect
to permutations of its arguments, it is called an \emph{Exchangeable
partition probability function (EPPF)}, and the corresponding
partition $\Pi_{n}$ an \emph{exchangeable random partition}.

 \noindent For exchangeable Gibbs partitions, the EPPF is, for
$\alpha\in(-\infty,1],$

\begin{equation}
q_{\alpha,V}(n_{1},\ldots,n_{k})={V_{n,k}\prod_{j=1}^{k}}(1-\alpha)_{(n_{j}-1)},
\label{geppf}
\end{equation}
 with $(V_{n,k})$ defined as in (\ref{v}). This can be obtained by
 repeated application of (\ref{gn})-(\ref{go})-(\ref{v}).

\noindent A minimal sufficient statistic for an exchangeable
$\Pi_{n}$ is given, because of the symmetry of its EPPF, by its
\emph{unordered frequencies} (i.e. the count of how many frequencies
in $\Pi_{n}$ are equal to 1, ..., to $n$), whose distribution is
given by their (unordered) sampling formula:

\begin{equation}
\tilde{\mu}(\mathbf{n})={|\mathbf{n}|\choose\mathbf{n}}\frac{1}{\prod_{1}^{n}b_{i}!}q_{\mu}(\mathbf{n}),
\label{unord}
\end{equation}
where
$${|\mathbf{n}|\choose\mathbf{n}}=\frac{|\mathbf{n}|!}{\prod_{j=1}^{k}n_{j}!}$$
\noindent and $b_{i}$ is the number of $n_{j}$'s in $\mathbf{n}$
equal to $i$ $(i=1,\ldots,n)$.

\noindent It is easy to see that for a Ewens' partition (whose EPPF
is (\ref{geppf}) with $\alpha=0$ and $V$ given by (\ref{ew})),
formula (\ref{unord}) returns the celebrated \emph{Ewens' sampling
formula}.

\noindent The distribution of the age-ordered frequencies

\begin{equation}
\bar{\mu}(\mathbf{n})={|\mathbf{n}|\choose\mathbf{n}}a(\mathbf{n})q_{\mu}(\mathbf{n}),
\label{ord}
\end{equation}
\ \\
differs from (\ref{unord}) only by a counting factor, where
\begin{equation}
a(\mathbf{n})=\prod_{j=1}^{k}\frac{n_{j}}{n-\sum_{i=1}^{j-1}n_{i}}
\label{sbp} \end{equation} is the distribution of the size-biased
permutation of $\mathbf{n}$.

 \noindent If  $\Pi=(\Pi_{n})$ is a (partially) exchangeable partition with PEPF $q_{\mu}$, then the vector $n^{-1}(n_{1},n_{2},\ldots)$
of the relative  frequencies, in age-order, of $\Pi_{n}$, converges
a.s. to a random point $P=(P_{1},P_2,\ldots)\in\Delta$ with
distribution $\mu$: thus the integrand in (\ref{eppf}) has the
interpretation as the conditional PEPF of a partially exchangeable
random partition, given its limit age-ordered frequencies
$(p_{1},p_2,\ldots)$. If $q_{\mu}$ is an EPPF, then the measure
$d\mu$ is \emph{invariant under size-biased permutation}.

\ \\

\noindent The notion of PEPF gives a generalized version of Hoppe's
urn scheme, i.e. a predictive distribution for (the age-ordered
frequencies of) $\Pi_{n+1}$, given (those of) $\Pi_{n}$.  In an urn
of Hoppe's type there are colored balls and a black ball. Every time
we draw a black ball, we return it in the urn with a new ball of a
new distinct color. Otherwise, we add in the urn a ball of the same
color as the ball just drawn.

  \noindent Pitman's extended urn scheme works as follows.
  Let $q$ be a PEPF, and assume that
initially in the urn there is only the black ball. Label with $j$
the $j$-th distinct color appearing in the sample. After $n\geq 1$
samples, suppose we have put in the urn
$\mb{n}=(n_{1},\ldots,n_{k})$ balls of colors $1,\ldots,k$,
respectively, with colors labeled by their order of appearance. The
probability that the next ball is of color $j$ is

\begin{equation}
\mathbb{P}(\mathbf{n}+\mathbf{e}_{j}|\mathbf{n})=\frac{q(\mathbf{n}+\mathbf{e}_{j})}{q(\mathbf{n})}\mbb{I}(j\leq
k)+\left(1-\sum_{j=1}^{k}\frac{q(\mathbf{n}+\mathbf{e}_{j})}{q(\mathbf{n})}\right)\mathbb{I}(j=k+1),
\label{sss}
\end{equation}
where $\mb{e}_j=(\delta_{ij}:i=1,\ldots,k)$ and $\delta_{xy}$ is the
Kronecker delta. The event $(j=k+1)$, in the last term of the
right-hand side of (\ref{sss}), corresponds to a new distinct color
being added to the urn.

\noindent The predictive distribution  of a Gibbs partition is
obtained from its EPPF by substituting (\ref{geppf}) into
(\ref{sss}):

\begin{equation}
\mathbb{P}(\mathbf{n}+\mathbf{e}_{j}|\mathbf{n})=\frac{n_{j}-\alpha}{n-\alpha
k}\left(1-\frac{V_{n+1,k+1}}{V_{n,k}}\right)\mbb{I}(j\leq
k)+\frac{V_{n+1,k+1}}{V_{n,k}}\mbb{I}(j=k+1), \label{hoppe3}
\end{equation}
which gives back our definition (\ref{gn})-(\ref{go}) of an EGP.

\noindent The use of an urn scheme of the form (\ref{gn})-(\ref{go})
in Population Genetics  is due to Hoppe \cite{H84} in the context of
Ewens' partitions (infinitely-many-alleles model), for which the
connection between order of appearance in a sample and age-order of
alleles is shown by \cite{D86}. In \cite{DK99} an extended version
of Hoppe's approach is suggested for more complicated, still
exchangeable population models (where e.g. mutation can be
recurrent). Outside Population Genetics, the use of
(\ref{gn})-(\ref{go}) for generating trees leading to Pitman's
two-parameter GEM frequencies, can be found in the literature of
random recursive trees (see e.g.  \cite{DGM06}). Urn schemes of the
form (\ref{sss}) are a most natural tool to express one's \emph{a
priori} opinions in a Bayesian Statistical context, as pointed out
by \cite{Z} and \cite{P96b}. Examples of recent applications of
Exchangeable Gibbs partitions in Bayesian nonparametric Statistics
are in \cite{LPW07}, \cite{IJ03}. The connection between (not
necessarily infinite) Gibbs partitions and coagulation-fragmentation
processes is explored by \cite{BP06} (see also \cite{B06} and
reference therein).

\subsection{Distribution of record indices in partially exchangeable
partitions.}

Let $\Pi$ be a partially exchangeable random partition. Since, for
every $n$, its age-ordered frequencies $\mb{n}=(n_{1},\ldots,n_{k})$
are a sufficient statistic for $\Pi_{n}$, then all realizations
$\pi_{n}$ with the same $\mb{n}$ \emph{and} the same record indices
$\mb{i}_{n}=1<i_{2}<\ldots<i_{k}\leq n$ must have equal probability.
To evaluate the joint probability of the pair $(\mb{n},\mb{i}_{n})$,
we only need to replace $a(\mb{n})$ in (\ref{ord}) by an appropriate
counting factor.
 This is equal to the
number of arrangements of $n$ balls, labelled from $1$ to $n$, in
$k$ boxes with the constraint that exactly $n_{j}$ balls fall in the
same box as the ball $i_{j}$. Such a number was shown by \cite{GL}
to be equal to

\begin{equation*}
{|\mathbf{n}|\choose\mathbf{n}}a(\mathbf{n},\mb{i}_{n}) \label{jc}
\end{equation*}
where
\begin{equation*}
{|\mathbf{n}|\choose\mathbf{n}}=\frac{n!}{\prod_{j=1}^{k}n_{j}!}\ \
\text{and}\ \
a(\mathbf{n},\mb{i}_{n})=\frac{\prod_{j=1}^{k}{{S_{j}-i_{j}}\choose{n_{j}-1}}}{{|\mathbf{n}|\choose\mathbf{n}}}
\normalcolor
\end{equation*} with $S_{j}:=\sum_{i=1}^{j}n_{i}$.
Thus, if $\Pi=(\Pi_{n})$ is a partially exchangeable random
partition with PEPF $q_{n}$, then the joint probability of
age-ordered frequencies and record indices is
\begin{equation}
\bar{\mu}(\mathbf{n},\mathbf{i}_{n})={|\mathbf{n}|\choose\mathbf{n}}a(\mathbf{n},\mb{i}_{n})q_{\mu}(\mathbf{n}).
\label{w}\end{equation}  The distribution of the record indices can
be easily derived by marginalizing:
\begin{equation}
\bar{\mu}(\mathbf{i}_{n})=\sum_{\mb{n}\in
B(\mb{i}_{n})}{|\mathbf{n}|\choose\mathbf{n}}a(\mathbf{n},\mb{i}_{n})q_{\mu}(\mathbf{n}).
\label{mi}\end{equation} where
$$B_{n}(\mb{i}_{n})=\{(n_{1},\ldots,n_{k}):\sum_{i=1}^{k}n_i=n;\ S_{j-1}\geq i_{j}-1, j=1,\ldots,k\}$$ is the set of all
possible $\mb{n}$ compatible with $\mb{i}$. In \cite{GL} such a
formula is derived for the particular case of Ewens' partitions. For
general random partitions see also \cite{N06}, section 2.

\noindent Notice that, for every $\mb{n}$ such that $|\mb{n}|=n$,
$$a(\mb{n})=\sum_{\mb{i}_{n}\in C(\mb{n})}a(\mb{n},\mb{i}_{n})=\prod_{j=1}^{k}\frac{n_{j}}{n-\sum_{i=1}^{j-1}n_{i}},$$
where $$C(\mb{n})=\{(1< i_{2}<\ldots<i_{k}\leq n):k\leq n,\
i_{j}\leq S_{j-1}+1\}$$ is the set of all possible $\mb{i}_{n}$
compatible with $\mb{n}$. Then the marginal distribution of the
age-ordered frequencies (\ref{ord}) is recovered by summing
(\ref{w}) over $C(\mb{n})$.

 \noindent This observation incidentally links a classical combinatorial result to partially exchangeable
random partitions.

\begin{proposition}
\label{geom} Let $\Pi=(\Pi_{n})$ be a partially exchangeable
partition with PEPF $q_{\mu}$. \\(i) Given the frequencies
$\mathbf{n}=(n_{1},\ldots,n_{k})$ in age-order, the probability that
the least elements of the classes of $\Pi_{n}$ are
$\mb{i}_{n}=(i_{1},\ldots,i_{k})$, does not depend on $q_{\mu}$ and
is given by
\begin{equation}
\mathbb{P}(\mb{i}_{n}|\mathbf{n})=\frac{1}{n!}\prod_{j=1}^{k}\frac{(S_{j}-i_{j})!}{(S_{j-1}-i_{j}+1)!}(n-S_{j-1}).
\label{condi}
\end{equation}
(ii) Let $W_{j}=\lim_{n\rightarrow\infty}S_{j}/n$. Conditional on
 $\{W_j:j=1,2,\ldots\}$, the waiting times
$$T_{j}=i_{j}-i_{j-1}-1\hspace{2cm}(j=2,3,\ldots)$$ are independent geometric random variables, each with parameter $(1-W_{j-1})$,
respectively.
\end{proposition}

\begin{proof}
  Part (i) can be obtained by a manipulation of a standard result on uniform random permutations of $[n]$. Part two can be proved by using a representation theorem due to Pitman (\cite{P95}, Theorem 8).
We prefer to give a direct proof of both parts to make clear their
connection.
 Simply notice that, for every $\mb{n}$ and $\mb{i}$, the right-hand side of
 (\ref{condi}) is equal to ${a(\mb{n},\mb{i}_{n})}/{a(\mb{n})}$.
 Then, for every $\mb{n}$,
$$\sum_{\mb{i}_{n}}\mathbb{P}(\mb{i}_{n}|\mb{n})=\sum_{\mb{i}_{n}\in C(\mb{n})}\frac{a(\mb{n},\mb{i}_{n})}{a(\mb{n})}=1$$
and
$$\sum_{\mb{n}}\sum_{C(\mb{n})}
\mathbb{P}(\mb{i}_{n}|\mb{n})\bar{\mu}(\mb{n})=1,$$ where
$\bar{\mu}(\mb{n})$ is as in (\ref{ord}), hence
$\mathbb{P}(\mb{i}_{n}|\mb{n})$ is a regular conditional probability
and (i) is proved. Now, consider the set
$$C_{[i_l,l]}(\mb{n}):=\{i_{l}<i_{l+1}<\ldots<i_{k}\leq n:i_{j}\leq
S_{j-1}+1,\ j=l+1,\ldots,k\}.$$ Also define, for $j=1,\ldots,k-l+1,$
$n_{j}^{*}:=S_{j}^{*}-S_{j-1}^{*}$ with
$S_{j}^{*}:=S_{j+l-1}-i_{l}+1,$ and $i_j^{*}:=i_{j+l-1}-i_l+1.$ Then
$C_{[i_l,l]}(\mb{n})=C(\mb{n}^{*})$  so that, for a fixed $l\leq k$,
the conditional probability of $i_2,\ldots,i_l$, given
$\mb{n}=(n_1,\ldots,n_k)$, is
\begin{eqnarray}
\mbb{P}(i_{2},\ldots,i_l|\mb{n})&=&\frac{1}{n!}\sum_{C_{[i_l,l]}(\mb{n})}\prod_{j=1}^{k}\frac{(S_{j}-i_{j})!}{(S_{j-1}-i_{j}+1)!}(n-S_{j-1})\nonumber\\
&=&\frac{(n-i_{l})!}{n!}\left[\prod_{j=1}^{l-1}\frac{(S_{j}-i_{j})!}{(S_{j-1}-i_j+1)!}(n-S_{j-1})\right]\frac{(n-S_{l-1})}{(S_{l-1}-i_l+1)!}\nonumber\\
&&\ \ \ \times
\sum_{C(\mb{n}^{*})}\frac{a(\mb{n}^{*},\mb{i}_{n^{*}}^{*})}{a(\mb{n}^{*})},\label{gpr1}
\end{eqnarray}
where $n^{*}=n-i_l+1$. The sum in (\ref{gpr1}) is 1; multiply and
divide the remaining part by $[S_{l}^{l-1}(S_{l}-i_l)!]/(S_l-1)!.$
The probability can therefore be rewritten as
\begin{equation*}
\mbb{P}(i_{2},\ldots,i_l|\mb{n})=\frac{(S_{l}-1)_{[i_l-1]}}{(n-1)_{[i_l-1]}}\left(\frac{S_{l}}{n}\right)^{-(l-1)}
\prod_{j=1}^{l}\left(1-\frac{S_{j-1}}{n}\right)\left[\frac{S_{l}^{l-1}}{(S_{l}-1)!}\prod_{j=1}^{l}\frac{(S_{j}-i_j)!}{(S_{j-1}-i_j+1)!}\right],
\end{equation*}
where $a_{[r]}=a(a-1)\cdots(a-r+1)$ is the falling factorial.
 Now, define
$$({{W}_{j}}:j=1,2,\ldots)=\lim_{n\rightarrow\infty}({{S}_{j}}/{n}:j=1,2,\ldots);$$
then, for $l$ fixed,
\begin{eqnarray*}
\lim_{n\rightarrow\infty}\mbb{P}(i_{2},\ldots,i_l|\mb{n})&=&W_{l}^{i_l-l}\prod_{j=1}^{l}(1-W_{j-1})\left[\prod_{j=2}^{l}\left(\frac{W_{j-1}}{W_{l}}\right)^{i_{j}-i_{j-1}-1}\right]\nonumber\\
&=&\prod_{j=2}^{l}W_{j-1}^{i_{j}-i_{j-1}-1}(1-W_{j-1}),
\end{eqnarray*}
 which is the distribution of $k-1$ independent
geometric random variables, each with parameter $(1-W_{j})$, and the
proof is complete.
\end{proof}

\bigskip

\noindent By combining the definition (\ref{eppf}) of PEPF and
Proposition \ref{geom}, one recovers an identity due to Nacu
(\cite{N06}, (7)).

\bigskip

\begin{corollary}
\text{(\cite{N06}, Proposition 5)} For every sequence
$1=i_1<i_{2}<\ldots<i_k+1=n$, and every point
$x=(x_1,x_2,\ldots)\in\Delta$,
\begin{equation}
\prod_{j=1}^{k-1}w_{j}^{i_{j+1}-i_{j}-1}=\sum_{B({\mb{i}})}{|\mathbf{n}|\choose\mathbf{n}}a(\mb{i},\mb{n})\prod_{j=1}^{k}x_{j}^{n_{j}-1}
\label{nacu}
\end{equation}
where $w_{j}=\sum_{i=1}^{j}x_{i}$ $(j=1,2,\ldots)$.
\end{corollary}

\begin{proof}
Multiply both sizes by $\prod_{1}^{k}(1-w_{j-1})$: by Proposition
\ref{geom} and (\ref{eppf}), formula (\ref{nacu}) is just the
equality (\ref{mi}) with the choice $d\mu=\delta_{x}$.
\end{proof}

\bigskip

\section{Age-ordered frequencies conditional on the record indices in Exchangeable
Gibbs partitions}

\subsection{Conditional distribution of sample
frequencies.}\label{*}

\noindent From now on we will focus only on EGP$(\alpha, V)$. We
have seen that the conditional distribution of the record indices,
given the age-ordered frequencies of a partially exchangeable random
partition, is purely combinatorial as it does not depend on its
PEPF. We will now find the conditional distribution of the
age-ordered frequencies $\mb{n}$ given the record indices, i.e. the
step (i) of the plan outlined in the introduction. We show that such
a distribution does not depend on the parameter $V$, which in fact
affects only the marginal distribution of $\mb{i}_{n}$, as explained
in the following Lemma.

\bigskip

\begin{lemma}
\label{margi} Let $\Pi=(\Pi_{n})$ be an EGP$(\alpha,V)$, for some
$\alpha\in(-\infty,1)$ and $V=(V_{n,k}:k\leq n=1,2,\ldots)$.
 For each $n$, the probability that the record indices in $\Pi_{n}$ are
$\mb{i}_{n}=(i_{1},\ldots,i_{k})$ is
\begin{equation}
\bar{\mu}_{\alpha,V}(\mb{i}_{n})=\psi_{\alpha,n,k}^{-1}(\mb{i}_{n})V_{n,k}.
\label{ri}
\end{equation}
where
\begin{equation}
\psi_{\alpha,n,k}(\mb{i}_{n})=\frac{\Gamma(1-\alpha)}{\Gamma(n-\alpha
k)}\prod_{2}^{k}\frac{\Gamma(i_{j}-j\alpha)}{\Gamma(i_{j}-j\alpha-(1-\alpha))};
\label{psi}
\end{equation}
Then the sequence $i_1,i_2,\ldots$ forms a non-homogeneous Markov
chain starting at $i_{1}=1$ and with transition probabilities given
by
\begin{equation}
P_{j}(i_{j+1}|i_{j})=(i_{j}-\alpha
j)_{(i_{j+1}-i_{j}-1)}\frac{V_{i_{j+1},j+1}}{V_{i_{j},j}},\hspace{1.5cm}j\geq
1. \label{itrp}
\end{equation}
\end{lemma}

\begin{proof}The proof can be carried out by using the urn scheme (\ref{hoppe3}).
For every $n$, let $K_{n}$ be the number of distinct colors which
appeared before the $n+1$-th ball was picked. From (\ref{hoppe3}),
the sequence $(K_{n}:n\geq 1)$ starts from $K_{1}=1$ and obeys, for
every $n$, the prediction rule:

\begin{equation}
\mathbb{P}(K_{n+1}=k_{n+1}|K_{n}=k)=\left(1-\frac{V_{n+1,k}}{V_{n,k}}\right)\mb{I}(k_{n+1}=k)+\frac{V_{n+1,k+1}}{V_{n,k}}\mb{I}(k_{n+1}=k+1).
\label{khoppe}
\end{equation}
By definition, $K_{n}$ jumps at points $1<i_{2}<\ldots$, due to the
equivalence
$$\{K_{n+1}=K_{n}+1|K_{n}=k\}=\{i_{k+1}=n+1\}.$$
Therefore, from (\ref{khoppe}),

\begin{eqnarray}
\bar{\mu}_{\alpha,V}(\mb{i}_{n})&=&\prod_{j=1}^{k}\frac{V_{i_{j},j}}{V_{i_{j-1},j-1}}\prod_{l\notin\mathbf{i}_{n}}\left(l-1-K_{l-1}\alpha\right)
\frac{V_{l,K_{l-1}}}{V_{l-1,K_{l-1}}}
\nonumber\\
&=&V_{n,k}\prod_{l\notin\mathbf{i}_{n}}\left(l-1-K_{l-1}\alpha\right).\label{g2i1}\end{eqnarray}
The last product in (\ref{g2i1}) is equal to

\begin{eqnarray}
\prod_{l\notin\mb{i}_{n}}\left(l-K_{l-1}\alpha-1\right)&=&\frac{\Gamma(n-\alpha
k)}{\Gamma(1-\alpha)}\prod_{j=2}^{k}\frac{\Gamma(i_{j}-j\alpha-(1-\alpha))}{\Gamma(i_{j}-j\alpha)}.
\nonumber\\
&=&\psi_{\alpha,n,k}^{-1}(\mb{i}_{n}).\label{im}
\end{eqnarray}
 and this proves (\ref{ri}). The second part of the Lemma (i.e. the
 transition probabilities (\ref{itrp})) follow immediately just by replacing, in
 (\ref{ri}), $n$ with $i_{k}$, for every $k$, to show that
 $$\bar{\mu}_{\alpha,V}(\mb{i}_{i_{k}})=\prod_{j=1}^{k-1}P_{j}(i_{j+1}|i_{j}),$$
 for $P_{j}$ satisfying (\ref{itrp}) for every $j$.
 \end{proof}

\bigskip

\noindent The distribution of the age-ordered frequencies in
an EGP $\Pi_{n}$, conditional on the record indices, can be easily
obtained from Lemma \ref{margi} and (\ref{w}).

\bigskip

\begin{proposition} \label{sample} Let $\Pi=(\Pi_{n})$ be an EGP$(\alpha,V)$, for
some $\alpha\in(-\infty,1)$ and $V=(V_{n,k}:k\leq n=1,2,\ldots)$.
 For each $n$, the conditional distribution of the sample
frequencies $\mb{n}$ in age-order, given the vector $\mb{i}_{n}$ of
indices, is independent of $V$ and is equal to

\begin{equation}
\bar{\mu}_{\alpha}(\mb{n}|\mb{i}_{n})=\psi_{\alpha,n,k}(\mb{i}_{n})\left(\prod_{j=1}^{k}{{S_{j}-i_{j}}\choose{n_{j}-1}}(1-\alpha)_{(n_{j}-1)}\right).
\label{gcond}
\end{equation}

\end{proposition}
\bigskip

\noindent \begin{remark} Notice that, as $\alpha\rightarrow 0$,
formula (\ref{gcond}) reduces to that for Ewens' partitions, proved
in \cite{GL}:

$$\bar{\mu}_{0}(\mb{n}|\mb{i}_{n})={\prod_{2}^{k}(i_{j}-1)}\prod_{j=1}^{k}\frac{({S_{j}-i_{j}})!}{(S_{j-1}+i_{j}-1)!}.$$
\end{remark}
\bigskip

\begin{proof}

\noindent Recall that the probability of a pair
$(\mathbf{n},\mathbf{i}_n)$ is given by

\begin{equation}
\bar{\mu}_{\alpha,V}(\mathbf{n},\mb{i}_{n})={|\mathbf{n}|\choose\mathbf{n}}a(\mathbf{n},\mb{i}_{n})q_{\alpha,V}(\mathbf{n}).
\label{2glsf}
\end{equation}
 Now it is easy to derive the conditional distribution of a
configuration given a sequence $\mathbf{i}_n$, as:

\begin{eqnarray}
\bar{\mu}_{\alpha}(\mathbf{n}|\mb{i}_{n})&=&\frac{\bar{\mu}_{\alpha,V}(\mathbf{n},\mb{i}_{n})}{\bar{\mu}_{\alpha,V}(\mb{i}_{n})}\nonumber\\
&=&{|\mathbf{n}|\choose\mathbf{n}}\left(\frac{\prod_{j=1}^{k}{{S_{j}-i_{j}}\choose{n_{j}-1}}(1-\alpha)_{(n_{j}-1)}}
{{|\mathbf{n}|\choose\mathbf{n}}}\right)\frac{V_{n,k}}{\bar{\mu}_{\alpha,V}(\mathbf{i})}\nonumber\\
&=&\psi_{\alpha,n,k}(\mb{i}_{n})\left(\prod_{j=1}^{k}{{S_{j}-i_{j}}\choose{n_{j}-1}}(1-\alpha)_{(n_{j}-1)}\right),
\label{2cond}
\end{eqnarray}
and the proof is complete.
\end{proof}

\bigskip

\subsection{The distribution of the limit frequencies given the
record indices.} \label{**}

 We now have all elements to derive a representation for the limit
relative frequencies in age-order, conditional on the limit sequence
of record indices $\mb{i}=(i_{1}<i_{2}<\ldots)$ generated by an
EGP$(\alpha,V)$.

\begin{proposition}
\label{cgem} Let $\Pi=(\Pi_{n})_{n\geq 1}$ be an Exchangeable Gibbs
Partition with index $\alpha>-\infty$ for some $V$. Let
$\mb{i}=(i_{1}<i_{2}<\ldots)$ be its limit sequence of record
indices and $X_{1},X_{2},\ldots$ be the age-ordered limit
frequencies as $n\rightarrow\infty$.

\noindent A regular conditional distribution of $X_{1},X_{2},\ldots$
given the record indices is given by
\begin{equation}
X_{j}\overset{d}{=}\xi_{j-1}\prod_{m=j}^{\infty}(1-\xi_{m}),\hspace{2cm}j\geq
1, \label{gglcondlim}
\end{equation}
a.s., where $\xi_{0}\equiv 1$ and, for $j\geq 1$, $\xi_{j}$ is a
Beta random variable in $[0,1]$ with parameters
$(1-\alpha,i_{j+1}-j\alpha-1)$.
\end{proposition}
 \bigskip
\begin{remark}\label{gnedrem}
Proposition \ref{cgem} is a statement about a regular conditional
distribution. The question about the existence of a limit
conditional distribution of $X|\mb{i}$ as a function of
$\mb{i}=\lim_{n}\mb{i}_{n}$ has different answer according to the
choice of $\alpha$, as a consequence of the limit behavior of $K_n$,
the number of blocks of an EGP $\Pi_n$, as recalled in the
introduction. For $\alpha<0$, $\mb{i}$ is almost surely a finite
sequence; for nonnegative $\alpha$, the length $l$ of $\mb{i}$ will
be a.s. either $l\sim s \log n$ (for $\alpha=0$) or $l\sim
sn^{\alpha}$ (for $\alpha>0$), for some $s\in[0,\infty]$. The
infinite product representation (\ref{gglcondlim}) still holds in
any case if we adopt the convention $i_k\equiv\infty$ for every
$k>K_{\infty}$ where $K_{\infty}:=\lim_{n\rightarrow\infty}K_n$.
\end{remark}

\bigskip
\begin{proof}

 The form (\ref{gcond}) of the conditional density
$\mu_{\alpha}(\mb{n}|\mb{i}_{n})$ implies
\begin{equation}
\sum_{|\mb{n}|=n}\prod_{j=1}^{k}{{S_{j}-i_{j}}\choose{n_{j}-1}}(1-\alpha)_{(n_{j}-1)}=\psi_{\alpha,n,k}^{-1}(\mb{i}_{n}).
\label{d1}\end{equation} For some $r<k$, let $a_{2},\ldots,a_{r}$ be
positive integers and set $a_{1}=0$ and $a_{r+1}=\ldots=a_{k}=0$.
Define $\mb{i}_{n}'=(i_{1}',\ldots,i_{k}')$ where
$$i_{j}'=i_{j}+\sum_{1}^{j}a_{i}\hspace{1.5 cm}(j=1,\ldots,k).$$
Now take the sum (\ref{d1}) with $\mb{i}_{n}$ replaced by
$\mb{i}_{n}'$, and multiply it by $\psi_{\alpha,n,k}(\mb{i}_{n})$.
We obtain
\begin{equation}
\frac{\psi_{\alpha,n,k}(\mb{i}_{n})}{\psi_{\alpha,n,k}(\mb{i}_{n}')}=\mathbb{E}\left(\prod_{j=1}^{k}\frac{(S_{j}-i_{j}')!}{(S_{j}-i_{j})!}\frac{(S_{j-1}-i_{j}+1)!}{(S_{j-1}-i_{j}'+1)!}\right)
\label{d2}\end{equation} where the expectation is taken with respect
to $\bar{\mu}_{\alpha}(\cdot |\mb{i}_{n})$. The left hand side of
(\ref{d2}) is
\begin{eqnarray}
\frac{\psi_{\alpha,n,k}(\mb{i}_{n})}{\psi_{\alpha,n,k}(\mb{i}_{n}')}&=&
\prod_{j=2}^{k}\frac{[i_{j}-j\alpha-(1-\alpha)]_{(\sum_{i=1}^{j}a_{i})}}{(i_{j}-j\alpha)_{(\sum_{i=1}^{j}a_{i})}}\nonumber\\
&=&\prod_{j=1}^{k-1}\mathbb{E}((1-\xi_{j})^{\sum_{i=1}^{j}a_{i}}),
\label{d3}
\end{eqnarray}
where $\xi_{1},\ldots,\xi_{k-1}$ are independent Beta random
variables, each with parameters $(1-\alpha,i_{j+1}-j\alpha-1)$.

\noindent Let $b_{j}=\sum_{i=1}^{j}a_{i}$. On the right hand side of
(\ref{d2}), $S_{0}=0, S_{k}=n$, so the product is equal to

\begin{equation}
\prod_{j=1}^{k}\left(\frac{S_{j-1}}{S_{j}}\right)^{b_{j}}\left[\prod_{l=0}^{b_{j}-1}\left(\frac{1-\frac{i_{j}+l}{S_{j}}}{1-\frac{i_{j}-1+l}{S_{j-1}}}\right)\right]
\label{d4}\end{equation} Since $a_{j}=0$ for $j=1$ and $j>r$, as
$k,n\rightarrow\infty$ the product inside square brackets
converges to 1 so the limit of (\ref{d4}) is
$$\prod_{j=1}^{r-1}W_{j}^{a_{j+1}}$$
where $W_{j}=\lim_{n\rightarrow\infty}S_{j}/n$. Hence from
(\ref{d2}) it follows that in the limit

$$\mathbb{E}\left(\prod_{j=1}^{r-1}W_{j}^{a_{j}}\right)=\prod_{j=2}^{\infty}E((1-\xi_{j})^{\sum_{i=1}^{j}a_{i}})$$
which gives the limit distribution of the cumulative sums:

\begin{equation*}
W_{j}\overset{d}{=} \prod_{i=j}^{\infty}(1-\xi_{i}),
\hspace{1.5cm}j=1,2,\ldots\end{equation*} But
\begin{eqnarray*}
X_{j}&=&W_{j}-W_{j-1}\\
&\overset{d}{=}&\prod_{i=j}^{\infty}(1-\xi_{j})-\prod_{i=j-1}^{\infty}(1-\xi_{j})\\
&=&\xi_{j-1}\prod_{i=j}^{\infty}(1-\xi_{i}),\hspace{1.5cm}j=1,2,\ldots
\end{eqnarray*}
and the proof is complete.
\end{proof}

\subsection{Conditional Gibbs frequencies, Neutral distributions and invariance under size-biased permutation.}

\noindent Proposition \ref{cgem} says that, conditional on all the
record indices $i_{1},i_{2},\ldots$, the sequence of relative
increments of an EGP$(\alpha,V)$

\begin{equation}
\xi=\left(\frac{X_{2}}{W_{2}},\frac{X_{3}}{W_{3}},\ldots\right).
\label{incr} \end{equation}
 is a sequence of independent
coordinates. In fact, such a process can be interpreted as the
negative, time-reversed version of a so-called \textit{Beta-Stacy}
process, a particular class of random discrete distributions,
introduced in the context of Bayesian Nonparametric Statistics
 as a useful tool to make inference for right-censored
data (see \cite{WM97}, \cite{WM99} for a modern account).\\
It is possible to show that such an independence property of the
$\xi$ sequence (conditional on the indices) actually characterizes
the family of EGP partitions. To make clear such a statement we
recall a concept of \emph{neutrality} for random $[0,1]$-valued
sequences, introduced by Connor and Mosimann \cite{CM} in 1962 and
refined in 1974 by Doksum \cite{Dk74} in the context of
nonparametric inference and, more recently, by Walker and Muliere
\cite{WM99}.

\begin{definition}

\label{ntr}  Let $k$ be any fixed positive integer (non necessarily finite). \\
\ \\
(i) Let $P=(P_{1},P_{2},\ldots,P_k)$ be a random point in
$[0,1]^{k}$ such that $\sum_{i=1}^{k}P_{i}\leq 1$ and, for every
$j=1,\ldots,k-1$ denote $F_{j}=\sum_{i=1}^{j}P_{i}$. P is called a
\emph{Neutral to the Right (NTR)} sequence if the vector
$B_{1},B_{2},\ldots,B_{k-1}$ of relative increments
\begin{equation*}
B_{j}=\frac{P_{j}}{1-F_{j-1}}
\hspace{2cm}j=1,\ldots,k-1\end{equation*} is a sequence of
independent random variables in ${[0,1]}$. \\
Let $(\alpha,\beta)$ be a point in $[0,\infty]^{k-1}$. A NTR vector
such that every increment $B_j$ is a Beta $(\alpha_j,\beta_j)$
$(j=1,\ldots,k-1)$, is called a \emph{Beta-Stacy} distribution with
parameter
$(\alpha,\beta)$. \\

 \noindent (ii) A \emph{Neutral to the left (NTL)}
vector $P=(P_{1},P_{2},\ldots,P_k)$ is a vector such that
$P^{*}:=(P_{k},P_{k-1},\ldots,P_{1})$ is NTR. \\
A \emph{Left-Beta-Stacy} distribution with parameter
$(\alpha,\beta)$ is a NTL vector $P$ such that $P^{*}$ is
Beta-Stacy
$(\alpha^{*},\beta^{*})=((\alpha_{k-1},\beta_{k-1}),\ldots,(\alpha_{1},\beta_{1}))$.
\end{definition}

\bigskip
\noindent A known result due to \cite{P96} is that the only class
of exchangeable random partitions whose limit age-ordered
frequencies are (unconditionally) a NTR distribution, is Pitman's
two-parameter family, i.e. the EGP$(\alpha,V)$ with
$V$-coefficients given by (\ref{p2p}). In this case, the
age-ordered frequencies follow the so-called two-parameter GEM
distribution, a special case of Beta-Stacy distribution with each
$B_{j}$ being a Beta$(1-\alpha,\theta+j\alpha)$ random
variable.\\
The age-ordered frequencies of all other Gibbs partitions are not
NTR; on the other side, Proposition \ref{cgem} shows that,
conditional on the record indices $i_1,i_2,\ldots$, and on $W_{k}$
they are all NTL distributions. For a fixed $k$ set
\begin{equation*}
 Y_{j}=\frac{X_{k-j+1}}{W_{k}},\hspace{1.5cm}1\leq j\leq
 k.\end{equation*}
Then
$${1-F_{j}}=\frac{W_{k-j}}{W_{k}}$$ and
\begin{equation}
\xi_{k-j}=\frac{Y_{j}}{1-F_{j-1}}.
 \label{xibs}
 \end{equation} By construction, the sequence $Y_{1},\ldots,Y_{k}$ is a Beta-Stacy sequence with
 parameters $\alpha_{k,j}=1-\alpha$ and
 $\beta_{k,j}=i_{k-j+1}-(k-j+1)\alpha-(1-\alpha)$ $(j=1,\ldots,k)$.
The property of (conditional) left-neutrality is maintained as
$k\rightarrow\infty$
 (just condition on $W_{k_{\infty}}=1$ where
 $K_{\infty}=\lim_{n\rightarrow\infty}K_{n}$ ).\\

 \noindent The
 following proposition is a converse of Proposition \ref{cgem}.

 \begin{proposition}
 \label{ntlp}
 Let $X=(X_{1},X_{2},\ldots)\in\Delta$ be the age-ordered frequencies of an
 infinite exchangeable random partition $\Pi$ of $\mbb{N}$. Assume,
 conditionally on the record indices of $\Pi$, $X$ is a NTL
 sequence.
 Then $\Pi$ is an exchangeable Gibbs partition for some
 parameters $(\alpha,V)$.
 \end{proposition}

\begin{proof}
 \noindent
 The frequencies of an exchangeable random partition of $\mbb{N}$ are in age-order if and only if their distribution is \emph{invariant
 under size-biased permutation} (ISBP, see \cite{DJ}, \cite{P96}).
 To prove the proposition, we
 combine two known results: the first is a characterization of  ISBP distributions; the second is a
 characterization of the Dirichlet distribution in terms of NTR
 processes. We recall such results in two lemmas.
 \bigskip
 \begin{lemma}\label{isbp}\  \emph{Invariance under size-biased permutation (\cite{P96},
 Theorem 4).} Let $X$ be a random point of $[0,1]^{\infty}$ such that
 $|P|\leq 1$ almost surely with respect to a probability measure $d\mu$. For every $k$, let $\mu_{k}$ denote the distribution
 of $X_{1},\ldots,X_{k},$ and $G_{k}$ the measure on $[0,1]^{k}$,
 absolutely continuous with respect to $\mu_{k}$ with density
 $$\frac{dG_{k}}{d\mu_{k}}(x_{1},\ldots,x_{k})=\prod_{i=1}^{k-1}(1-w_{j})$$
 where $w_j=\sum_{i=1}^{j}x_{i}, \ j=1,2,\ldots$ \\
 $X$ is invariant under size-biased permutation if and only if
 $G_{k}$ is symmetric with respect to permutations of the
 coordinates in $\mbb{R}^{k}$.
 \end{lemma}

\noindent Let $X$ be the frequencies of an exchangeable partition
$\Pi$, and denote with $\mbb{P}_{\mu}$ the marginal law of the
record indices of $\Pi$. Consider the measure $G_{k}$ of Lemma
\ref{isbp}. By Proposition \ref{geom} (ii), for every $k$
\begin{eqnarray}
G_{k}(dx_{1}\times\cdots\times dx_{k})&=&
\mu_{k}(dx_{1}\times\cdots\times
dx_{k})\prod_{i=1}^{k-1}(1-w_{j})\nonumber\\
&&=\mu_{k}(dx_{1}\times\cdots\times
dx_{k})\mbb{P}(i_1=1,i_{2}=2,\ldots,i_{k}=k\ | x_{1},\ldots,x_{k})\nonumber\\
&&=\mu_{k}(dx_{1}\times\ldots\times
dx_{k}|i_{k}=k)\mbb{P}_{\mu}(i_{k}=k)\label{isbp2}
\end{eqnarray}
An equivalent characterization of ISBP measure is:

\begin{corollary}\label{isbpcor}
The law of $X$ is invariant under size-biased permutation if and
only if, for every $k$, there is a version of the conditional
distribution
\begin{equation*}
\mu_{k}(dx_{1}\times\ldots\times dx_{k}|i_{k}=k)
\end{equation*}
which is invariant under permutations of coordinates in
$\mbb{R}^{k}.$
\end{corollary}
\bigskip

\noindent The other result we recall is about Dirichlet
distributions.

\bigskip

\begin{lemma}
\label{ntrdir} \ \emph{Dirichlet and neutrality (\cite{BW},
Theorem 7).} Let $P$ be a random $k$-dimensional vector with
positive components such that their sum equals 1. If $P$ is NTR
and $P_{n}$ does not depend on
$(1-X_{k})^{-1}(P_{1},\ldots,P_{k})$. Then $X$ has the Dirichlet
distribution.
\end{lemma}

\bigskip

\noindent Now we have all elements to prove Proposition
\ref{ntlp}. Let $\mu(\cdot | i_{1},i_{2}\ldots)$ be the
distribution of a NTL vector such that the distribution of
$\xi_j:=X_{j+1}/W_{j+1}$ has marginal law $d\gamma_j$ for
$j=1,2,\ldots$. For every $k$, given $i_1,\ldots,i_k$, the vector
$(X_{2}/W_{2},\ldots,X_{k}/W_{k})$ is conditionally independent of
$W_{k}$ and
\begin{equation}
\mu_k(dx_{1}\times\cdots\times dx_k |
i_{k}=k)=\left(\prod_{j=1}^{k-1}\gamma_j(d\xi_{j})\right)\zeta_k(dw_{k})\nonumber\
\end{equation}
where $\zeta_k$ is the conditional law of $W_k$ given $i_{k}=k$.\\
For $X$ to be ISBP, corollary \ref{isbpcor} implies that the
product
$$\prod_{j=1}^{k-1}\gamma_j(d\xi_{j})$$
must be a symmetric function of $x_1,\ldots,x_k.$ Then, for every
$k$, the vector $(X_{1}/W_{k},\ldots,X_{k}/W_{k})$ is both NTL and
NTR, which implies in particular that $X_{k}/W_k$ is independent
of $W_{k-1}^{-1}(X_1,\ldots,X_{k-1})$. Therefore, by Lemma
\ref{ntrdir} and symmetry,
$(\frac{X_1}{W_k},\ldots,\frac{X_{k}}{W_{k}})$ is, conditionally
on $W_{k}$ and $\{i_k=k\}$, a symmetric Dirichlet distribution,
with parameter, say, $1-\alpha>0.$ By (\ref{isbp2}), the EPPF
corresponding to $d\mu$ is equal to
\begin{equation}
\mbb{E}\left[\prod_{j=1}^{k}X_{j}^{n_j-1}(1-W_{j-1})\right]=\mbb{P}_{\mu}(i_{k}=k)\mbb{E}\left[\prod_{j=1}^{k}X_{j}^{n_j-1}|i_{k}=k\right].\label{xeq1}\end{equation}
By the NTL assumption, we can write
\begin{eqnarray}
\prod_{j=1}^{k}X_{j}^{n_j-1}&=&\prod_{j=2}^{k}\left(\frac{X_{j}}{W_{j}}\right)^{n_j-1}\left(\frac{W_{j}}{W_{k}}\right)^{n_{j}-1}W_{k}^{n-k},\nonumber\\
&\overset{d}{=}&
\left(\prod_{j=2}^{k}\xi_{j}^{n_{j+1}-1}(1-\xi_{j})^{S_{j}-j}\right)W_{k}^{n-k},\label{xeq2}
\end{eqnarray}
where $S_{j}=\sum_{i=1}^{j}n_{i}$ $(j=1,\ldots,k)$. The last
equality is due to
$$\prod_{j=1}^{k}\frac{W_{j}}{W_{k}}=\prod_{j=1}^{k-1}\left(\frac{W_{j}}{W_{j+1}}\right)^{j}.$$
Now, set
$$V_{n,k}=\frac{\mbb{P}_{\mu}(i_{k}=k)\mbb{E}(W^{n-k}_{k}|i_{k}=k)}{[k(1-\alpha)]_{(n-k)}}.$$ Equality (\ref{xeq2}) implies
\begin{eqnarray*}
\mbb{E}\left[\prod_{j=1}^{k}X_{j}^{n_j-1}(1-W_{j-1})\right]&=&\mbb{P}_{\mu}(i_{k}=k)\int\left(\int\prod_{j=1}^{k-1}\xi_{j}^{n_{j+1}-1}(1-\xi_{j})^{S_j-j}\gamma_{j}(\xi_{j})d\xi_{j}\right)w_{k}^{n-k}\zeta_{k}(w_{k})dw_{k}\nonumber\\
&=&\mbb{P}_{\mu}(i_{k}=k)\mbb{E}(W^{n-k}_{k}|i_{k}=k)\prod_{j=1}^{k-1}\frac{(1-\alpha)_{(n_{j+1}-1)}[j(1-\alpha)]_{(S_j-j)}}{[(j+1)(1-\alpha)]_{(S_{j+1}-(j+1))}}\nonumber\\
&=&V_{n,k}\prod_{j=1}^{k}(1-\alpha)_{(n_{j}-1)},
\end{eqnarray*}
which completes the proof.
\end{proof}
\bigskip

\section{Age-ordered frequencies conditional on a single record index.}
 \label{***}
A representation for the Mellin transform of the $m$-th
age-ordered cumulative frequencies $W_{m}$, conditional on $i_{m}$
alone $(m=1,2,\ldots)$ can be derived by using Proposition
\ref{cgem}. \noindent   We first point out a characterization for
the moments of $W_m$, stated in the following Lemma.

\begin{lemma}
\label{momprp} Let $X_{1},X_{2},\ldots$ be the limit age-ordered
frequencies generated by a Gibbs partition with parameters
$\alpha,V$. For every $m=1,2,\ldots$ and nonnegative integer $n$
\begin{equation}
\mathbb{E}(W_{m}^{n}|i_{m})=(i_{m}-\alpha
m)_{(n)}\frac{V_{i_{m}+n,m}}{V_{i_{m},m}}\label{wmom}
\end{equation}
and
\begin{equation}
\mbb{E}(X_m^n|i_m)=(1-\alpha)_{(n)}\frac{V_{i_{m}+n,m}}{V_{i_{m},m}}.
\label{fmom}
\end{equation}
\end{lemma}

\begin{proof}
Let $\Pi$ be  an EGP $(\alpha,V)$ and denote $Y_{j}:j=1,2,\ldots$
the sequence of indicators $\{0,1\}$ such that $Y_{j}=1$ if $j$ is a
record index of $\Pi$. Then $Y_{1}=1$ and, for every $l,m\leq l$,
\begin{eqnarray*}
\mbb{P}(Y_{l+1}=0\ |\ \sum_{i=1}^{l}Y_{i}=m)&=&(l-\alpha
m)\frac{V_{l+1,m}}{V_{l,m}}\nonumber\\
&=&1-\mbb{P}(Y_{l+1}=1\ |\ \sum_{i=1}^{l}Y_{i}=m).
\end{eqnarray*}
 By proposition \ref{geom} and formula (\ref{geppf}), given the cumulative frequencies $W= W_{1},W_{2},\ldots$,
$$\mbb{P}(Y_{l+1}=0|\sum_{i=1}^{l}Y_{i}=m, W)=1-\mbb{P}(Y_{l+1}=1\ |\ \sum_{i=1}^{l}Y_{i}=m, W)=W_{m}.$$
(see also \cite{P95}, Theorem 6). Obviously this also implies that,
conditional on $W$, the random sequence $K_{l}:=\sum_{i=1}^{l}Y_i,\
(l=1,2,\ldots)$ is Markov , so we can write, for every $l,m$

\begin{equation}
\mbb{P}(Y_{l+1}=0\ |\ K_l=m,\ K_l-1=m-e,W)=\mbb{P}(Y_{l+1}=0\ |\
K_l=m,W)=W_{m},\hspace{1cm}e=0,1.\label{mw}
\end{equation}

\noindent Hence

\begin{eqnarray}\mbb{E}(W_{m}|i_{m})&=&\mbb{E}\left[\mbb{P}(Y_{i_m+1}=0\ |\
K_{i_{m}}=m,K_{i_{m}-1}=m-1,W)\right]\nonumber\\
&=&\mbb{E}\left[\mbb{P}(Y_{i_{m}+1}=0\ |\ K_{i_{m}}=m,W)\right]\nonumber\\
 &=&(i_m-\alpha m)\frac{V_{i_m+1,m}}{V_{i_m,m}}\label{oh},
 \end{eqnarray}
which proves the proposition for $m=1$. The Markov property of $K_n$
and (\ref{mw}) also
 lead, for every $m$, to
\begin{eqnarray*}\mbb{E}(W^{n}_{m} |i_{j})&=&\mbb{E}\left[\mbb{P}(Y_{i_m+1}=\ldots=Y_{n}=0\ |\ K_{i_{m}}=m,K_{i_{m}-1}=m-1,W)\right]\nonumber\\
&=&\prod_{i=1}^{n}\mbb{E}\left[\mbb{P}(Y_{i_m+i}=0\ |\ K_{i_{m}+i-1}=m,W)\right]\nonumber\\
 &=&(i_m-\alpha m)_{(n)}\frac{V_{i_m+n,m}}{V_{i_m,m}},
 \end{eqnarray*}
where the last equality is obtained as an $n$-fold iteration of
(\ref{oh}).


\noindent The second part of the Lemma (formula (\ref{fmom}))
follows from Proposition \ref{cgem}:
\begin{eqnarray*}
\mbb{E}(X_m^n|i_m)&=&\mbb{E}(\xi_{m-1}^{n}|i_m)\mbb{E}(W_{m}^{n}|i_{m})\nonumber\\
&=&\frac{(1-\alpha)_{(n)}}{(i_m-\alpha
m)_{(n)}}\mbb{E}(W_{m}^{n}|i_{m})
\end{eqnarray*}
which combined with (\ref{wmom}) completes the proof.
\end{proof}

\noindent Given the coefficients $\{V_{n,k}\}$, analogous formulas
 to (\ref{wmom}) and (\ref{fmom}) can be obtained to describe the
conditional Mellin transforms of $W_{m}$ and $X_{m}$
(respectively), in terms of the Mellin transform of the
size-biased pick $X_1$.

\bigskip

\begin{proposition}\label{melprp}
Let $X_{1},X_{2},\ldots$ be the limit age-ordered frequencies
generated by a Gibbs partition with parameters $\alpha,V$. For every
$m=1,2,\ldots$ and $\phi\geq 0$
\begin{equation}
\mathbb{E}(W_{m}^\phi|i_{m})=(i_{m}-\alpha
m)_{(\phi)}\frac{V_{i_{m},m}[\phi]}{V_{i_{m},m}}\label{wmel}
\end{equation}
and
\begin{equation}
\mbb{E}(X_m^\phi|i_m)=(1-\alpha)_{(\phi)}\frac{V_{i_{m},m}[\phi]}{V_{i_{m},m}},
\label{fmel}
\end{equation}
for a sequence of functions
$(V_{n,k}[\cdot]:\mbb{R}\rightarrow\mbb{R};k,n=1,2,\ldots)$,
uniquely determined by $V_{n,k}\equiv V_{n,k}[0]$,
 such that, for
every $\phi\geq 0,$
\begin{eqnarray}
&&V_{1,1}[\phi]=\frac{\mbb{E}(X_{1}^{\phi})}{(1-\alpha)}_{(\phi)};\label{v11}\\
&& \nonumber\\
&& V_{n,k}[\phi]=(n+\phi-\alpha
k)V_{n+1,k}[\phi]+V_{n+1,k+1}[\phi]\label{vnk},\hspace{1cm}n,k=1,2,\ldots;\\
&& \nonumber\\
&&
V_{n,k}[\phi+1]=V_{n+1,k}[\phi]\hspace{1cm}n,k=1,2,\ldots.\label{addv}
\end{eqnarray}
\end{proposition}

\begin{remark}
 To
complete the representation given in Proposition \ref{melprp},
notice that, for every $\alpha$, the distribution of $X_{1}$ (the
so-called \emph{structural distribution}) is known for the extreme
points of the Gibbs$(\alpha,V)$ family.  \\
In particular,
for $\alpha\leq 0$ $X_1$ has a Beta$(1-\alpha,\theta+\alpha)$
density ($\theta>0$), where $\theta=m|\alpha|$ for some integer
$m$ when $\alpha<0$ (see e.g. \cite{P96b}). In this case,
$$V_{1,1}[\phi](\theta)=\frac{1}{(1-\alpha)_{(\phi)}}\left(\frac{(1-\alpha)_{(\phi)}}{(\theta+1)_{(\phi)}}\right)=\frac{1}{(\theta+1)_{(\phi)}}.$$
When $\alpha>0$, we saw in the introduction that every extreme
point in the Gibbs family is a Poisson-Kingman $(\alpha,s)$
partition for some $s>0$; in this case the density of $X_1$ is

$$f_{1,\alpha}(x|s)=\frac{\alpha s^{-1}(x)^{-\alpha}}{\Gamma(1-\alpha)}\frac{f_{\alpha}((1-x)s^{1/\alpha})}{f_{\alpha}(s^{1/\alpha})}\hspace{2cm}0<x<1,$$
for an $\alpha$-stable density $f_\alpha$ (\cite{P92}, (57)),
leading to
$$V_{1,1}[\phi](s)=\alpha
s^{\frac{\phi-1}{\alpha}}G_\alpha(\phi-\alpha-1,s^{-1/\alpha})$$
where $G_\alpha$ is as in (\ref{galpha}). \\
Thus, for every $\alpha$, the structural distribution of a
Gibbs$(\alpha,V)$ partition, which defines $V_{1,1}[\phi]$ in
(\ref{v11}), can be obtained as mixture of the corresponding
extreme structural distributions.

\noindent
\end{remark}

\begin{proof}
Note that, for $\phi=0,1,2,\ldots$, the proposition holds by Lemma
\ref{momprp} with $V_{n,m}[\phi]\equiv V_{n+\phi,m}$. For general
$\phi\geq 0$ observe that, for every $m,n\in\mbb{N}$,
\begin{equation}
\mbb{E}(W_{m}^{\phi+n-m}|i_{m}=m)=\mbb{E}(W_{m}^{\phi}|i_{m}=n)\mbb{E}(W_{m}^{n-m}|i_{m}=m).
\label{idmel}
\end{equation}
To see this, consider the random sequences $Y_{n},K_n$ defined in
the proof of Lemma \ref{momprp}. By (\ref{mw}),
\begin{eqnarray*}
\mbb{E}(W_{m}^{\phi+n-m}|i_{m}=m)&=&\mbb{E}\left[W_{m}^{\phi}\mbb{P}(K_{n}=m|K_{m}=m,W)
\ |\ K_{m}=m\right] \nonumber\\
&=& \mbb{E}\left[W_{m}^{\phi}\mbb{P}(K_{n}=m|\ W) \ |\
K_{m}=m\right]
\nonumber\\
&=& \mbb{E}\left[W_{m}^{\phi} \ |\
K_{n}=m,K_{m}=m\right]\mbb{E}\left[\mbb{P}(K_{n}=m|\ W) \ |\
K_{m}=m\right]\nonumber\\
&=& \mbb{E}\left[W_{m}^{\phi} \ |\
K_{n}=m\right]\mbb{E}\left[W_m^{n-m} \ |\ K_{m}=m\right]\nonumber\\
&=&\mbb{E}\left[W_{m}^{\phi} \ |\ i_m=n\right]\mbb{E}\left[W_m^{n-m}
\ |\ K_{m}=m\right]
\end{eqnarray*}
where the last two equalities follow from (\ref{mw}), the Markov
property of $K_{n}$ and the exchangeability of the $Y's$.

\noindent From Lemma \ref{momprp}, we can rewrite (\ref{idmel}) as
\begin{equation}
\mbb{E}(W_{m}^{\phi}|i_{m}=n)=\frac{\mbb{E}(W_{m}^{\phi+n-m}|i_{m}=m)}{[m(1-\alpha)]_{(n-m)}}\frac{V_{m,m}}{V_{n,m}}.
\label{idmel2}
\end{equation}
Now define
\begin{equation}
M_{m}(\phi)=\mbb{E}(W_{m}^{\phi}|i_{m}=m),\hspace{1cm}\phi\geq
0,m=1,2,\ldots \label{vmmphi} .\end{equation}
 and
\begin{equation}
V_{n,m}[\phi]=\frac{V_{m,m}M_{m}(\phi+n-m)}{[m(1-\alpha)]_{(n+\phi-m)}}\hspace{1cm}\phi\geq
0,n,m=1,2,\ldots. \label{vnmphi}
\end{equation}
Notice that, with such a definition, Lemma 1 implies that
$V_{n,m}[0]=V_{n,m}$. Moreover, $V_{n,m}[\phi+1]=V_{n+1,m}[\phi]$
so (\ref{addv}) is satisfied; then (\ref{idmel2}) reads
\begin{equation}
\mathbb{E}(W_{m}^\phi|i_{m})=(i_{m}-\alpha
m)_{(\phi)}\frac{V_{i_{m},m}[\phi]}{V_{i_{m},m}},
\end{equation}
that is: (\ref{v11}),(\ref{wmel}) and  are satisfied.
\\
Now it only remains to prove that such choice of $V_{n,m}[\phi]$
obeys the recursion (\ref{vnk}) for every $n,m,\phi.$ By the same
arguments leading to (\ref{idmel}),

\begin{eqnarray}
M_{m}(\phi)-M_{m}(\phi+1)&=&\mbb{E}[W_m^{\phi}(1-W_m)|K_m=m]\nonumber\\
&=&\mbb{E}\left[W_{m}^{\phi}\big(1-\mbb{P}(K_{m+1}=m|W) \big)\ |\
K_{m}=m\right]\nonumber\\
&=&\mbb{E}\left[W_{m}^{\phi} \ |\ K_{m+1}=m+1,
K_{m}=m\right]\mbb{E}\left[1-\mbb{P}(K_{m+1}=m|W) \ |\
K_{m}=m\right]\nonumber\\
&=&\mbb{E}\left[(1-\xi_{m})^{\phi} \ |\ i_{m+1}=m+1
\right]\mbb{E}\left[W_{m+1}^{\phi} \ |\ K_{m+1}=m+1
\right]\frac{V_{m+1,m+1}}{V_{m,m}}. \label{wd}
\end{eqnarray}
The last equality is a consequence of Proposition \ref{cgem}, for
which, conditional on $K_{m+1}=m+1$ (which is equivalent to
$\{i_{m+1}=m+1\}$), $W_{m}=(1-\xi_{m})W_{m+1}$ for $\xi_m$
(independent of $W_{m+1}$) having a Beta$(1-\alpha, m(1-\alpha))$
distribution. Therefore (\ref{wd}) can be rewritten as
\begin{equation}
M_{m}(\phi)-M_{m}(\phi+1)=\frac{[m(1-\alpha)]_{(\phi)}}{[(m+1)(1-\alpha)]_{(\phi)}}M_{m+1}(\phi)\frac{V_{m+1,m+1}}{V_{m,m}},
\label{mkrec}\end{equation} and (\ref{vnk}) follows from
(\ref{mkrec}),(\ref{addv}), after some simple algebra, by
comparing the definition (\ref{vnmphi}) of $V_{n,k}[\phi]$, and
the recursion (\ref{v}) for the $V$-coefficients of an
EGP$(\alpha,V)$. In particular, (\ref{mkrec}) shows that the
functions $V_{n,k}[\phi]$ are uniquely determined by
$V=(V_{n,k})$.\\
The equality (\ref{fmel}) can now be proved in the same way as the
moment formula (\ref{fmom}).
\end{proof}

\bigskip

\subsection{Convolution structure.} \label{***i}

An alternative representation of $-\log X_{m}$ can be given in terms
of an infinite sum of random variables, as done by \cite{GL} for
Kingman's model (Ewens' partition). Formula (42) in \cite{GL} shows
that, conditional on $i_{m}$
\begin{equation}
-\log(X_{j})=\sum_{i=2}^{i_{j}}Y_{i}+\sum_{m=i_{j}+1}^{\infty}Z_{m},
\label{gllog}
\end{equation}
where $\{Y_{i}\}$ is a sequence of independent exponential random
variables such that the $i$-th component has rate $i-1$, and
$\{Z_{m}\}$ is a sequence of mutually independent random variables
with density continuous everywhere in the non-negative reals except
one atom at zero:

\begin{equation}
f_{i}(w)=\frac{1}{1+\rho_{i}}\mb{I}(w=0)+\frac{\rho_{i}}{1-\rho_{i}}(i-1)e^{-(i-1)w}\mb{I}(w>0),
\label{glcont}
\end{equation}
where $ \rho_{i}=\frac{\theta}{i-1}$ (\cite{GL}, p. 173).

\noindent With the help of (\ref{gglcondlim}) it is possible to show
that a similar representation holds true for general Gibbs
partitions, but that the collection of the $Y$'s is split in two subsequences. It turns out that
 the independence property of the $Z$'s actually
characterizes Kingman's model: for general Gibbs partitions with
parameters $\alpha,V$, the $Z$'s are stochastically
linked by the sequence $K=(K_{n}:n=1,2,\ldots)$ as defined in the
proof of Proposition \ref{sample}. Remember that $K$ is a Markov
process on $\mathbb{N}$ such that $K_{1}=1$ and for $j>1$

\begin{eqnarray}
p_{kj}=P(K_{j}=k+1|K_{j-1}=k)&=&1-P(K_{j}=k|K_{j-1}=k)\nonumber\\
&=&\frac{V_{j,k+1}}{V_{j-1,k}}.\label{pkl}
\end{eqnarray}

\begin{proposition}
Let $X=X_{1},X_{2},\ldots$ be the limit relative frequencies
generated by a Gibbs partition with parameters $\alpha,V$.
Conditionally on $i_{m}$,
\begin{equation}
-\log
X_{m}=\sum_{j=2}^{m}Y^{*}_{j}+\sum_{i=m+1}^{i_{m}}Y^{(m)}_{i}+\sum_{l=i_m+1}^{\infty}Z_{l}
\label{logprop}
\end{equation}
where: $\{Y^{*}_{j}\}$ are independent random variables, with
distribution, respectively,
\begin{equation}
f^{(\alpha)}_{j}(y)dy=\frac{e^{-y(j-1-\alpha
(j-1))}{(1-e^{-y})}^{-\alpha}}{B(j-1-\alpha
(j-1),1-\alpha)}\mathbb{I}(y>0)dy; \label{y1}
\end{equation}
$\{Y^{(m)}_{i}\}$'s are independent random variables, with
distribution, respectively,
\begin{equation}
f^{(\alpha,m)}_{i}(y)dy=\frac{e^{-y(i-1-\alpha m)}}{i-\alpha
m-1}\mathbb{I}(y>0)dy, \label{ym}
\end{equation}
and $\{Z_{l}\}$ are such that such that: (i) $Z_{l}$ is
conditionally independent of $Z_{l-1}$, given $K_{l-1}$; (ii)
conditionally on $\{K_{l-1}=k,K_{l}=v\}$, the density of $Z_{l}$
is

\begin{equation}
g^{\alpha,V}_{l,k,v}(z)dz=\delta_{kv}\delta_{0z}dz+(1-\delta_{kv})\frac{e^{-z(l-1-\alpha
k)}{(1-e^{-z})}^{-\alpha}}{B(l-1-\alpha
k,1-\alpha)}\mathbb{I}(z>0)dz. \label{w0}
\end{equation}

\end{proposition}

\begin{proof}
From Proposition (\ref{cgem}), and the properties of Beta
distributions,

\begin{eqnarray}
&&\ \mathbb{E}(e^{\phi\log
X_{m}}|i_{m})=\mathbb{E}((\xi_{m-1}^{\phi}\prod_{i=m}^{\infty}(1-\xi_{i})^{\phi})|i_{m})\nonumber\\
&&\hspace{2cm}=\frac{(1-\alpha)_{(\phi)}}{(i_{m}-\alpha
m)_{(\phi)}}\mathbb{E}\left(\prod_{k=m}^{\infty}\frac{(i_{k+1}-\alpha
k-1)_{(\phi)}}{(i_{k+1}-\alpha(k+1))_{(\phi)}}\ \Big | \
i_{m}\right). \label{plap1}
\end{eqnarray}
 The fraction on the left-hand side can be decomposed as:
\begin{eqnarray}
\frac{(1-\alpha)_{(\phi)}}{(i_{m}-\alpha m)_{(\phi)}}&=&\left(\prod_{j=2}^{m}\frac{(j-1-\alpha(j-1))_{(\phi)}}{(j-\alpha
j)_{(\phi)}}\right)\left(\prod_{i=m+1}^{i_{m}}\frac{(i-1-\alpha
m)_{(\phi)}}{(i-\alpha m)_{(\phi)}}\right)\nonumber\\
&=&\left(\prod_{j=2}^{m}\mathbb{E}(\wt{Y}^{*\phi}_{j})\right)\left(\prod_{i=m+1}^{i_{m}}\mathbb{E}(\wt{Y}^{(m)\phi}_{i})\right),
\label{2p}\end{eqnarray} where, for each $j$,
$$\wt{Y}^{*}_{j}\sim {\rm Beta}(j-1-\alpha(j-1),1-\alpha),$$ and for each $i$,
$$\wt{Y}^{(m)}_{i}\sim {\rm Beta}(j-1-\alpha m,1).$$
Now simply set ${Y}^{*}_{j}:=-\log\wt{Y}^{*}_{j}$ and
${Y}^{(m)}_{i}:=-\log \wt{Y}^{(m)}_{i}$ to see that
\begin{equation*}
\frac{(1-\alpha)_{(\phi)}}{(i_{m}-\alpha m)_{(\phi)}}=\mathbb{E}(e^{-\phi[\sum_{j=2}^{m}Y^{*}_{j}+\sum_{i=m+1}^{i_{m}}Y^{(m)}_{i}]}),
\end{equation*}
which provides the first two sums in the representation
(\ref{logprop}).

\noindent Define $\chi_{l}=K_{l}-K_{l-1}$ where $K_{l}$ is
described by (\ref{pkl}). Then
\begin{eqnarray*}
\mathbb{E}\left(\prod_{k=m}^{\infty}\frac{(i_{k+1}-\alpha
k-1)_{(\phi)}}{(i_{k+1}-\alpha(k+1))_{(\phi)}}|i_{m}\right)&=&
\mathbb{E}\left(\prod_{l=i_{m}+1}^{\infty}\frac{(l-\alpha
K_{l}-(1-\alpha))_{\phi\chi_{l}}}{(l-\alpha(k_{l}))_{\phi\chi_{l}}}|i_{m}\right)\\
&=&\mathbb{E}\left(\prod_{l=i_{m}+1}^{\infty}\mathbb{E}(\wt{Z}_{l}^{\phi\chi_{l}}|K_{l},K_{l-1})|i_{m}\right),
\end{eqnarray*}
where each $\wt{Z}^{l}$ is a Beta random variable with parameters
$(l-\alpha K_l-(1-\alpha),1-\alpha)$. Now define
$$Z_{l}=-\chi_{l}\log\wt{Z}_{l}.$$ Then $Z_{l}$ has the required
distribution (\ref{w0}) and the last expectation can be written as
$$\mathbb{E}\left(\prod_{l=i_{m}+1}^{\infty}\mathbb{E}(e^{-{\phi}{Z}_{l}}|K_{l},K_{l-1})|i_{m}\right)=\mathbb{E}\left(e^{-\phi\sum_{l=i_{m}+1}^{\infty}Z_{l}}|i_{m}\right)$$
which completes the proof.
\end{proof}

\noindent The characterization of Kingman's model follows
immediately.

\begin{corollary}
Let $X=X_{1},X_{2},\ldots$ be the limit relative frequencies
generated by a Gibbs partition with parameters $\alpha,V$. The
random variables involved in the three sums in (\ref{logprop}) are
all mutually independent if and only if $\alpha=0$ and $V$ is such
that, for every $n,k\leq n,$ $V_{n,k}=\theta^{k}/\theta_{(n)}$ for
some $\theta>0$.
\end{corollary}

\begin{proof}
We only have to verify the independence of the $Z$'s. In general,
$Z_l$ depends on $Z_{l-1}$ only through $K_{l-1}$ $(l>i_m)$. By
averaging over $K_l$, one has that the conditional density of
$W_l$ given $\{K_{l-1}=k\}$ is:
\begin{equation}
(1-p_{kl})\delta_{0z}(dz)+p_{kl}\frac{e^{-z(l-1-\alpha
k)}{(1-e^{-z})}^{-\alpha}}{B(l-1-\alpha
k,1-\alpha)}\mathbb{I}(z>0)dz. \label{uncz}
\end{equation}
But remember that
$$1-p_{kl}=(l-\alpha
k-1)\frac{V_{l,k}}{V_{l-1,k}}$$ does not depend on $k$ only if
$\alpha=0$ and
$$\frac{V_{l,k}}{V_{l-1,k}}=\frac{c_{l-1}}{c_{l}}\hspace{1.5 cm}k,l=1,2,\ldots$$ for some constant $c_{l}$.
This implies that the $V$-coefficients are of the form
$$V_{n,k}=c_{n}^{-1}V^{*}_{k};$$but, as mentioned in the
introduction, within the family of all EGPs,  only  Pitman's
two-parameter partitions have such a form. Since we just saw that
$\alpha$ is necessarily zero, then it must be
$$V_{n,k}=\frac{\theta^{k}}{\theta_{(n)}}$$ for some $\theta>0,$
which corresponds exactly to Kingman's model. In this case, the
density (\ref{uncz}) reduces to
\begin{equation}\frac{l-1}{l+\theta-1}\delta_{0z}(dz)+\frac{\theta}{l+\theta-1}\frac{e^{-z(l-1)}}{l-1}.\label{bobz}
\end{equation}

\noindent The converse is straightforward.
\end{proof}

\begin{remark}
Notice that, when $\alpha=0$, the $Y_j^*$'s and the $Y_j^{(m)}$ are
all exponential random variables, each with parameter $(j-1)$,
respectively. This, together with (\ref{bobz}), leads us back to
Griffiths and Lessard's representation (\ref{gllog})-(\ref{glcont}).
\end{remark}

\section{A representation for normalized age-ordered frequencies in an exchangeable Gibbs partition.}
\label{****}

 \noindent  In this section we provide a characterization of
  the density of the first $k$ (normalized) age-ordered frequencies, given $i_{k}$ and $W_{k}$,
  and an explicit formula for the marginal distribution of $i_{k}$. We give
a direct proof, obtained by comparison of the unconditional
distribution of $X_1,\ldots,X_k$, $(k=1,2,\ldots)$, in a general
Gibbs partition, with its analogue in Pitman's two-parameter
model. Such a comparison is naturally induced by proposition
\ref{cgem}, which says that, conditional on the record indices,
the distribution of the age-ordered frequencies is the same for
every Gibbs partition. Remember that the limit (unconditional)
age-ordered frequencies in such a family are described by the
two-parameter GEM distribution, for which
$$X_{j}\overset{d}{=}B_{j}\prod_{i=1}^{j-1}(1-B_{i})$$ for a sequence $B_{1},B_{2},\ldots$
of independent Beta random variables with parameters, respectively,
$(1-\alpha,\theta+j\alpha)$ (see e.g. \cite{P96b}).

\begin{proposition}
\label{bay}
Let $X_{1},X_{2},\ldots$ be the age-ordered frequencies of a
EGP$(\alpha,V)$ and, for every $k$ let $W_{k}=\sum_{i=1}^{k}X_{i}$.

 \noindent (i) Conditional on $W_{k}=w$ and on $i_{k}=k+i$, the law of the vector ${X}_{1},\ldots,{X}_{k}$
is
\begin{equation}
d\mu_{\alpha,V}({x}_{1},\ldots,{x}_{k}|w,k+i)=\frac{{V}_{k+i,k}}{V_{k+i-1,k-1}}w^{-(k-1)}\sum_{\mb{m}\in\mathbb{N}^{k-1}:|\mb{m}|=k+i-1}\bar{\mu}_{\alpha,V}(\mb{m})
\mathcal{D}_{\mb{(m-\alpha)}}(\frac{{x}_{1}}{w},\ldots,\frac{{x}_{k}}{w})d{x}_{1},\ldots,d{x}_{k-1}
\label{mixdir}
\end{equation}
where: $\mathcal{D}_{\mb{(m-\alpha)}}$ is the Dirichlet density with parameters $(m_{1}-\alpha,\ldots,m_{k-1}-\alpha,1-\alpha)$, and
 $\mu_{\alpha,V}(\mb{m})$ is the age-ordered sampling formula (\ref{ord}) with Gibbs' EPPF $q_{\alpha,V}$ as in (\ref{geppf}):
$$\bar{\mu}_{\alpha,V}(\mb{m})=\left(\frac{(k+i-1)!}{\prod_{j=1}^{k-1}(m_{j})!(k+i-1-\sum_{i=1}^{j-1}m
_{i})}\right){V_{k+i-1,k-1}\prod_{j=1}^{k-1}}(1-\alpha)_{(m_{j}-1)}.$$

\noindent(ii) The marginal distribution of $i_{k}$ is
\begin{equation}
\mathbb{P}(i_{k}=k+i)=V_{k+i,k}\frac{\alpha^{-(k-1)}}{(k+i-1)!}\sum_{j=0}^{k-1}\frac{(-1)^{j+k+i-1}}{j!(k-j-1)!}(\alpha j)_{[k+i-1]},
\label{imarg}
\end{equation}
where $a_{[n]}=a(a-1)\cdots(a-n+1)$.
\end{proposition}
 \bigskip

\begin{remark}
The density of $i_{k}$ can be expressed in terms of generalized
Stirling numbers ${n\brack k}_{\alpha}$, defined as the
coefficients of $x^{n}$ in
$$\frac{n!}{\alpha^{k}k!}(1-(1-x)^{\alpha})^{k}$$
( \cite{K},\cite{GP}). Formula (\ref{imarg}) can be re-expressed
as:
$$\mathbb{P}(i_{k})=V_{i_k,k}{i_k-1\brack
k-1}_{\alpha},$$ which makes clear the connection between the
distribution of $i_k$ and the distribution of $K_n$ recalled in
the introduction (formula (\ref{knpr})). In fact, (\ref{imarg})
can be deduced simply from (\ref{knpr}) by the Markov property of
the sequence $K_n$ as
\begin{eqnarray*}
\mathbb{P}(i_{k})&=&\mathbb{P}(K_{i_{k}}=k \ |\ K_{i_{k}-1}=k-1
)\ \mathbb{P}(K_{i_{k}-1}=k-1)\\
&=&\frac{V_{i_k,k}}{V_{i_{k-1},k-1}} V_{i_{k-1},k-1} {i_k-1\brack
k-1}_{\alpha}.
\end{eqnarray*}
However here we give a self-contained proof in order to show how
(\ref{imarg}) is implied by proposition \ref{cgem} through
(\ref{mixdir}).
\end{remark}

\begin{proof}
From (\ref{itr}), we know that
\begin{equation*}
P(i_{2},\ldots,i_{k})=V_{i_{k},k}\frac{\Gamma(i_{2}-\alpha-1) \cdots
\Gamma(i_{k}-\alpha(k-1)-1)}{\Gamma (1-\alpha)\Gamma(i_{2}-2\alpha)
\cdots \Gamma(i_{k-1}-(k-1)\alpha)};
\end{equation*}
By proposition \ref{cgem} and Lemma \ref{momprp},
\begin{eqnarray}
\mathbb{E}(\prod_{j=1}^{k}X_{j}^{n_{j}}|i_{1},\ldots,i_{k})&=&\mathbb{E}\left(\prod_{j=1}^{k-1}\xi_{j}^{n_{j+1}}(1-\xi_{j})^{S_{j}}\prod_{i=k}^{\infty}(1-\xi_{i})^{n}\ \Big | \ i_{1},\ldots,i_{k}\right)\nonumber\\
&=&
\prod_{j=1}^{k-1}\frac{(1-\alpha)_{(n_{j+1})}(i_{j+1}-j\alpha-1)_{(S_{j})}}{(i_{j+1}-(j+1)\alpha)_{(S_{j+1})}}\mathbb{E}(\prod_{i=k}^{\infty}(1-\xi_{i})^{n}|i_{k})\nonumber\\
&=&\left(\prod_{j=1}^{k-1}\frac{(1-\alpha)_{(n_{j+1})}(i_{j+1}-j\alpha-1)_{(S_{j})}}{(i_{j+1}-(j+1)\alpha)_{(S_{j+1})}}\right)\frac{\Gamma(i_{k}+n-\alpha
k)}{\Gamma(i_{k}-\alpha k)}\frac{V_{n+i_{k},k}}{V_{i_{k},k}},\nonumber\\
\label{open}
\end{eqnarray}
hence
\begin{eqnarray}
\mathbb{E}& ( & \prod_{j=1}^{k}X_{j}^{n_{j}}|1\equiv
i_{1},i_{2},\ldots,i_{k})P(i_{2},\ldots,i_{k})\nonumber\\
&=&(i_{k}-\alpha
k)_{(n)}V_{n+i_{k},k}\prod_{j=1}^{k-1}\frac{(1-\alpha)_{n_{j+1}}\Gamma(i_{j+1}-j\alpha-1+S_{j})\Gamma(i_{j+1}-\alpha(j+1))}{\Gamma(i_{j+1}-\alpha(j+1)+S_{j+1})\Gamma(i_{j}-\alpha
j)}\nonumber\\
&=&\prod_{j=1}^{k}(1-\alpha)_{(n_{j})}\left[V_{i_{k}+n,k}\prod_{j=1}^{k-1}(i_{j}-\alpha
j+S_{j})_{(i_{j+1}-i_{j}-1)}\right], \label{openv}\end{eqnarray}
where the last equality follows after multiplying and dividing all
terms by $(1-\alpha)_{n_{1}}$. Consequently, a moment formula for
general Gibbs partitions is of the form:
\begin{equation}
\mathbb{E}_{(\alpha,V)}\left(\prod_{j=1}^{k}X_{j}^{n_{j}}\right)=\prod_{j=1}^{k}(1-\alpha)_{(n_{j})}\sum_{1<i_2<\ldots<i_{k}}c_{(1,i_{2})}\cdots
c_{(i_{k-1},i_{k})}V_{n+i_{k},k}. \label{gmm}
\end{equation}
where, for $1<j\leq k-1$, $$c_{(i_{j},i_{j+1})}=(i_{j}-\alpha j +
S_{j})_{(i_{j+1}-i_{j}-1)}.$$ For fixed
 $i_{k}$ denote
$$\lambda_{i_{k}}=\sum_{1<i_{2}<\cdots<i_{k}}c_{(i_{1},i_{2})}\cdots
c_{(i_{k-1},i_{k})}.$$ Then (\ref{gmm}) reads
\begin{equation}
{\mathbb{E}_{(\alpha,V)}\left(\prod_{j=1}^{k}X_{j}^{n_{j}}\right)}={\prod_{j=1}^{k}(1-\alpha)_{(n_{j})}}\sum_{i=k}^{\infty}\lambda_{i}V_{n+i,k}.
\label{gmm2}\end{equation} For Pitman's two-parameter family, this
becomes
\begin{equation}{\mathbb{E}_{(\alpha,\theta)}\left(\prod_{j=1}^{k}X_{j}^{n_{j}}\right)}={\prod_{j=1}^{k}(1-\alpha)_{(n_{j})}}
\sum_{i=k}^{\infty}\lambda_{i}\frac{\prod_{j=1}^{k}(\theta+\alpha(j-1))}{\theta_{(n+i)}}.
\label{pitprv} \end{equation}
 The two-parameter GEM distribution implies that
\begin{eqnarray}
\mathbb{E}\left(\prod_{j=1}^{k}X_{j}^{n_{j}}\right)&=&\mathbb{E}\left(\prod_{j=1}^{k}B_{j}^{n_{j}}(1-B_{j})^{n-S_{j}}\right)\nonumber\\
&=&\prod_{j=1}^{k}\frac{(1-\alpha)_{n_{j}}(\theta+j\alpha)_{(n-S_{j})}}{(\theta+1+(j-1)\alpha)_{(n-S_{j-1})}}\nonumber\\
&=&\frac{1}{(\theta)_{(n)}}\prod_{j=1}^{k}\frac{(1-\alpha)_{(n_{j})}(\theta+\alpha(j-1))}{(\theta+\alpha(j-1)+n-S_{j-1})}.
\label{2mme}
\end{eqnarray}
therefore, from (\ref{pitprv}) we derive the identity:
\begin{equation}
\prod_{j=1}^{k}\frac{1}{(\theta+\alpha(j-1)+n-S_{j-1})}=
{\sum_{i=k}^{\infty}}\lambda_{i}\frac{\theta_{(n)}}{\theta_{(n+i)}},\hspace{1cm}\theta>0.
 \label{2andg}\end{equation}
For $\theta>0$ replace $\theta$ by $\theta-n$ in (\ref{2andg}), and denote $n^{*}_{j}=1-\alpha+n_{j}$. We now find an expansion of
the left-hand side of (\ref{2andg}) in terms of products of the type $\prod_{1}^{k}(n^{*}_{j})_{(m_{j})}$, for $m_1,\ldots,m_{k}\geq 0$.
The left-hand side of (\ref{2andg}) is now
\begin{eqnarray}
&& \prod_{j=1}^{k}\frac{1}{(\theta-n+\alpha(j-1)+n-S^{*}_{j-1})}=\prod_{j=1}^{k}\frac{1}{(\theta+j-1+S^{*}_{j-1})}\nonumber\\
&&\hspace{1cm}=\prod_{j=1}^{k}\int_{0}^{1}t_{j-1}^{j-1+\theta-1-S^{*}_{j-1}}dt_{j-1}\nonumber\\
&&\hspace{1cm}=\int \left(\prod_{j=1}^{k-1}t_{j}\right)^{\theta-1}(t_{1}\cdots t_{k-1})^{-n^{*}_{1}}
(t_{2}\cdots t_{k-1})^{-n^{*}_{2}}\cdots t_{k-1}^{-n^{*}_{k-1}}\prod_{j=1}^{k-1}t_{j}^{j}\ dt_{0}\cdots dt_{k-1},
\label{lhsexp}
\end{eqnarray}
where $S^{*}_{0}=0$ and $S^{*}_{j}=\sum_{i=1}^{j}n^{*}_{j},\ j=1,\ldots,k-1$.
Make the change of variable
$$u_{j}=1-\prod_{i=j}^{k-1}t_{i},\hspace{2cm}j=0,\ldots,k-1.$$
Then $0<u_{k-1}<\ldots<u_{0}<1$. The absolute value of the Jacobian is
\begin{eqnarray*}
\left|\frac{du}{dt}\right|&=&(t_{1}\cdots t_{k-1})\times (t_{2}\cdots t_{k-1})\times t_{k-1}\times 1\\
&=&\prod_{j=1}^{k-1}t_{j}^{j}.
\end{eqnarray*}
Thus (\ref{lhsexp}) transforms to
$$\int_{0<u_{k-1}<\ldots<u_{0}<1}(1-u_{0})^{\theta-1}\prod_{j=1}^{k-1}(1-u)^{-n^{*}_{j}}du_{1}\cdots du_{k-1} du_{0}.$$
Fix $u_{0}$ and consider
\begin{eqnarray*}
&&\int_{0<u_{k-1}<\ldots<u_{0}}\prod_{j=1}^{k-1}(1-u)^{m_{j}}du_{1}\cdots du_{k-1}\\
&&\hspace{1.5cm}=u_{0}^{k-1+\sum_{i=1}^{k-1}m_{i}}\frac{1}{m_{j-1}+1}\cdot\frac{1}{m_{j-1}+m_{j-2}+1}\cdot\frac{1}{m_{j-1}+\ldots+m_{1}+1} .
\end{eqnarray*}
The integral in (\ref{lhsexp}) is thus
\begin{equation}\sum_{m_{1},\ldots,m_{k-1}\geq 0}\frac{c(\mb{m})}
{\prod_{j=1}^{k-1}(1-\alpha)_{(m_{j})}}\prod_{j=1}^{k-1}{n^{*}_{j}}_{(m_{j})}
\int_{0}^{1}(1-u_{0})^{\theta-1}u_{0}^{k-1+\sum_{i=1}^{k-1}m_{i}}du_{0},
\label{4}
\end{equation}
where
\begin{equation}
c(\mb{m})=\prod_{j=1}^{k-1}\frac{1}{k-j+\sum_{i=j}^{k-1}m_{i}}\prod_{j=1}^{k-1}\frac{(1-\alpha)_{(m_{j})}}{m_{j}!}.
\label{cm}
\end{equation}
Now consider the right-hand side of (\ref{2andg}), again with $\theta$ replaced by $\theta-n$:
$$\sum_{i=0}^{\infty}\frac{\lambda_{k+i}}{(k+i-1)!}\int_{0}^{1}(1-x)^{\theta-1}x^{i+k-1}dx$$
and compare it with (\ref{4}) to obtain a representation for the $\lambda_{i}$'s in (\ref{gmm2}):
\begin{equation}\frac{\lambda_{k+i}}{(k+i-1)!}=\sum_{\mb{m}\in\mathbb{N}_{0}^{k}:|\mb{m}|=i}\frac{c(\mb{m})}{(1-\alpha)_{(m_{j})}}
\prod_{j=1}^{k-1}{n^{*}_{j}}_{(m_{j})},\label{lambdai}
\end{equation}
where $\mathbb{N}_{0}=\mathbb{N}\cup\{0\}$. Recall that $n^{*}_{j}=1-\alpha+n_{j}$ and consider the identity
$$\frac{(1-\alpha)_{(n_{j})}{n^{*}_{j}}_{(m_{j})}}{(1-\alpha)_{(m_{j})}}=(1-\alpha+m_{j})_{(n_{j})}.$$
From (\ref{wmom}) we also know that
$$V_{n+k+i,k}=\frac{V_{k+i,k}}{(k+i-\alpha k)_{(n)}}\mathbb{E}(W_{k}^{n}|i_{k}=k+i).$$
Thus (\ref{gmm2}) and (\ref{lambdai}) imply that
\begin{eqnarray}
\mathbb{E}_{(\alpha,V)}\left(\prod_{j=1}^{k}X_{j}^{n_{j}}\right)&=&\sum_{i=0}^{\infty}(k+i-1)!\mathbb{E}(W_{k}^{n}|i_{k}=k+i)V_{k+i,k}\nonumber\\
&&\times\sum_{\mb{m}\in\mathbb{N}_{0}^{k-1}:|\mb{m}|=i}c(\mb{m})
{\left[\frac{(1-\alpha)_{(n_{k})}\prod_{j=1}^{k-1}(1-\alpha+m_{j})_{(n_{j})}}{(k+i-\alpha
k)_{(n)}}\right]}. \label{temp}
\end{eqnarray}
Now, for every $\mb{m}\in \mathbb{N}_{0}^{k-1}$ such that $|\mb{m}|=i$, define $m{'}_{j}=m_{j}+1$; then we can rewrite
\begin{eqnarray*}
(k+i-1)!c(\mb{m})V_{k+i,k}&=&\frac{V_{k+i,k}}{V_{k+i-1,k-1}}\left(\frac{(k+i-1)!}{\prod_{j=1}^{k-1}m{'}_{j}!(k+i-1-\sum_{l=1}^{j-1}m{'}_{j})}\right)V_{k+i-1,k-1}\prod_{j=1}^{k-1}(1-\alpha)_{(m{'}_{j}-1)}\\
&=&\frac{V_{k+i,k}}{V_{k+i-1,k-1}}\mu_{\alpha,V}(\mb{m}{'}).
\end{eqnarray*}
Thus the right-hand side of (\ref{temp}) becomes
\begin{equation}
=\sum_{i=0}^{\infty}\frac{V_{k+i,k}}{V_{k+i-1,k-1}}\sum_{\mb{m}{'}\in\mathbb{N}^{k-1}:|\mb{m}{'}|=k+i-1}\mu_{\alpha,V}(\mb{m}{'})
{\left[\frac{(1-\alpha)_{(n_{k})}\prod_{j=1}^{k-1}(m_{j}{'}-\alpha)_{(n_{j})}}{(k+i-\alpha k)_{(n)}}\mathbb{E}(W_{k}^{n}|i_{k}=k+i)\right]}
\label{finali}
\end{equation}
The term between square brackets is the $n_{1},\ldots,n_{k}$-th moment of $k$ $[0,1]$-valued random variables $Y_{1}(\mb{m}'),\ldots,Y_{k}(\mb{m}')$
such that
$$\sum_{i=1}^{k}{Y_{i}(\mb{m}')}\overset{d}{=}(W_{k}|i_{k}=k+i)$$
and, conditional on $\sum_{i=1}^{k}{Y_{i}(\mb{m}')}$, the distribution of
$$\frac{Y_{1}(\mb{m}')}{\sum_{i=1}^{k}{Y_{i}(\mb{m}')}},\ldots,\frac{Y_{k}(\mb{m}')}{\sum_{i=1}^{k}{Y_{i}(\mb{m}')}}$$
is a Dirichlet distribution with parameters $(m_{1}{'}-\alpha,\ldots,m_{k-1}{'}-\alpha,1-\alpha)$, therefore (\ref{finali})
 completes the proof of part (i).

 \ \\

 \noindent To prove part (ii), we only have to notice that, from (\ref{finali}) it must follow that
$$\sum_{i=0}^{\infty}\frac{V_{k+i,k}}{V_{k+i-1,k-1}}\sum_{\mb{m}{'}\in\mathbb{N}^{k-1}:|\mb{m}{'}|=k+i-1}\mu_{\alpha,V}(\mb{m}{'})=1$$
hence a version of the marginal probability of $i_{k}$ is, for every $k$
\begin{equation}
\mathbb{P}(i_{k}=k+i)=\frac{V_{k+i,k}}{V_{k+i-1,k-1}}
\sum_{\mb{m}{'}\in\mathbb{N}^{k-1}:|\mb{m}{'}|=k+i-1}\mu_{\alpha,V}(\mb{m}{'}).
\label{pitemp}
\end{equation}
This can be also argued directly, simply by noting that, in an EGP$(\alpha,V)$,
$$\sum_{\mb{m}{'}\in\mathbb{N}^{k-1}:|\mb{m}{'}|=k+i-1}\mu_{\alpha,V}(\mb{m}{'})=\mathbb{P}(i_{k-1}\leq k+i-1,i_{k}>k+i-1)$$ and that
$$\frac{V_{k+i,k}}{V_{k+i-1,k-1}}=\mathbb{P}(i_{k}=k+i|i_{k-1}\leq k+i-1, i_{k}>k+i-1).$$
We want to find an expression for the inner sum of (\ref{pitemp}).
 If we reconsider the term $c(\mb{m})$ as in (\ref{cm}) (for $\mb{m}\in\mathbb{N}_{0}^{k-1}:|\mb{m}|=i$), we see that
$$\sum_{\mb{m}\in\mathbb{N}_{0}^{k-1}:|\mb{m}|=i}c(\mb{m})$$
is the coefficient of $\zeta^{i}$ in
\begin{eqnarray*}
\frac{1}{(k-1)!}\left[\int_{0}^{1}(1-u\zeta)^{\alpha-1}du\right]^{k-1}&=&\frac{1}{(k-1)!}(\frac{1}{\zeta\alpha})^{k-1}[1-(1-\zeta)^{\alpha}]^{k-1}\\
&=&(\zeta\alpha)^{-(k-1)}\sum_{j=0}^{k-1}\frac{(-1)^{k+i+j-1}}{j!(k-1-j)!}(1-\zeta)^{\alpha j}.
\end{eqnarray*}
Thus

\begin{equation}\sum_{\mb{m}\in\mathbb{N}_{0}^{k-1}:|\mb{m}|=i}
c(\mb{m})=\frac{\alpha^{-(k-1)}}{(k+i-1)!}\sum_{j=0}^{k-1}\frac{(-1)^{k+i+j-1}}{j!(k-1-j)!}(j\alpha)_{[k+i-1]}.
\label{expl}\end{equation}
Since
$$\sum_{\mb{m}{'}\in\mathbb{N}^{k-1}:|\mb{m}{'}|=k+i-1}\mu_{\alpha,V}(\mb{m}{'})=(k+i-1)!\ V_{k+i-1,k-1}\sum_{\mb{m}\in\mathbb{N}_{0}^{k-1}:|\mb{m}|=i}c(\mb{m})$$
then part (ii) is proved by comparison of (\ref{expl}) with (\ref{pitemp}).
\end{proof}

\bibliographystyle{abbrv}

\begin{thebibliography}{10}

 \bibitem{A85}
    D.J.~Aldous
\newblock {\em Exchangeability and related topics}, volume 1117 of {\em Lecture
  Notes in Mathematics}, pp. 1--198.
\newblock Springer-Verlag, Berlin, 1985.
\newblock Lecture notes from \'{E}cole d'\'{e}t\'{e} de Probabilit\'{e}s de
  Saint-Flour XIII - 1983.

\bibitem{B06}
J. Bertoin
\newblock {\em Random
fragmentation and coagulation processes}.
\newblock Cambridge
Studies in Advanced Mathematics, vol. 102. Cambridge University
Press, Cambridge 2006.

\bibitem{BP06}
N.~Berestycki and J.~Pitman.
\newblock Gibbs distributions for random partitions generated by a
  fragmentation process, 2005.
\newblock http://arXiv.org/abs/math/0512378.

\bibitem{CM}
R.J.~Connor and J.E. Mosimann.
\newblock Concepts of independence for proportions with a generalization of the {D}irichlet distribution.
\newblock {\em J. Am. Stat. Assoc.}, 64:194--206, 1969.

\bibitem{BW}
K.~Bobecka and J. Weso\l owski.
\newblock The Dirichlet distribution and process through neutralities.
\newblock {\em J. Theor. Probab.}, 20: 295–-308, 2007.

\bibitem{Dk74}
K.~Doksum.
\newblock Tailfree and neutral random probabilities and their posterior distributions.
\newblock {\em Ann. Probab.}, 2:183--201, 1974.

\bibitem{DGM06}
R.~Dong, C.~Goldschmidt, and J.~B. Martin.
\newblock Coagulation-fragmentation duality, {P}oisson-{D}irichlet
  distributions and random recursive trees, 2005.
\newblock To appear in \emph{Ann. Appl. Probab.},
  http://arxiv.org/abs/math.PR/0507591.

\bibitem{D86}
P.~Donnelly.
\newblock Partition structures, {P}\'olya urns, the {E}wens sampling formula,
  and the ages of alleles.
\newblock {\em Theoret. Population Biol.}, 30(2):271--288, 1986.

\bibitem{DJ}
P.~Donnelly and P.~Joyce.
\newblock Continuity and weak convergence of ranked and size-biased permutations on the infinite simplex.
\newblock {\em Stoc. Proc. Appl.},31(1):89--103, 1989.

\bibitem{DK99}
P.~Donnelly and T.~G. Kurtz.
\newblock Particle representations for measure-valued population models.
\newblock {\em Ann. Probab.}, 27(1):166--205, 1999.

\bibitem{G97}
A.~V.~Gnedin.
\newblock The representation of composition structures.
\newblock {\em Ann. Probab.}, 25(3):1437--1450, 1997.

\bibitem{G06}
A.~Gnedin.
\newblock Constrained exchangeable partitions, 2006.
\newblock To appear in \emph{Discr. Math. Comp. Sci.},
  http://arXiv:math/0608621v1 [math.PR]

\bibitem{GP}
A.~Gnedin and J.~Pitman.
\newblock Exchangeable {G}ibbs partitions and {S}tirling triangles.
\newblock {\em Zap. Nauchn. Sem. S.-Peterburg. Otdel. Mat. Inst. Steklov.
  (POMI)}, 325(Teor. Predst. Din. Sist. Komb. i Algoritm. Metody. 12):83--102,
  244--245, 2005.

\bibitem{GL}
R.~C. Griffiths and S.~Lessard.
\newblock Ewens' sampling formula and related formulae: Combinatorial proofs,
  extensions to variable population size and applications to ages of alleles.
\newblock {\em Theor. Popul. Biol.}, 68:167--177, 2005.

\bibitem{H84}
F.~M. Hoppe.
\newblock P\'olya-like urns and the {E}wens' sampling formula.
\newblock {\em J. Math. Biol.}, 20(1):91--94, 1984.

\bibitem{IJ03}
H. Ishwaran and L.F. James.
 \newblock Some further developments for stick-breaking priors: finite
              and infinite clustering and classification,
\newblock {\em Sankhy\=a. The Indian Journal of Statistics}, 65(3): 577--592, 2003.


\bibitem{K}
S.~V. Kerov.
\newblock Combinatorial examples in the theory of {AF}-algebras.
\newblock {\em Zap. Nauchn. Sem. Leningrad. Otdel. Mat. Inst. Steklov. (LOMI)},
  172(Differentsialnaya Geom. Gruppy Li i Mekh. Vol. 10):55--67, 169--170,
  1989.

\bibitem{K2}
S.~V. Kerov.
\newblock Subordinators and permutation actions with quasi-invariant
              measure.
\newblock {\em Zap. Nauchn. Sem. S.-Peterburg. Otdel. Mat. Inst. Steklov.
              (POMI)},
  223(Teor. Predstav. Din. Sistemy, Kombin. i Algoritm. Metody. I):181--218, 340,
  1995.

\bibitem{KT1}
S.~V. Kerov and N.V. Tsilevich.
\newblock A random subdivision of an interval generates virtual
              permutations with the {E}wens distribution.
\newblock {\em Zap. Nauchn. Sem. S.-Peterburg. Otdel. Mat. Inst. Steklov. (POMI)},
  223(Teor. Predstav. Din. Sistemy, Kombin. i Algoritm. Metody. I):162--180, 339--340,
  1995.

\bibitem{K82}
J.~F.~C. Kingman.
\newblock On the genealogy of large populations.
\newblock {\em J. Appl. Probab.}, (Special Vol. 19A):27--43, 1982.
\newblock Essays in statistical science.

\bibitem{LPW07}
A. Lijoi, I. Pr\"{u}nster, S.G. Walker.
\newblock Bayesian nonparametric estimators derived from conditional
Gibbs structures.\emph{Tech. Report, Universit\`{a} degli Studi di
Torino, 2007}.

\bibitem{N06}
S.~Nacu.
\newblock Increments of random partitions, 2004.
\newblock http://arXiv:math/0310091v2 [math.PR]

\bibitem{PPY}
M.~Perman, J.~Pitman and M.~Yor.
\newblock Size-biased sampling of {P}oisson point processes and excursions
\newblock {\em Probab. Theory Related Fields}, 92(1):21--39, 1992.

\bibitem{P95}
J.~Pitman.
\newblock Exchangeable and partially exchangeable random partitions.
\newblock {\em Probab. Theory Related Fields}, 102(2):145--158, 1995.

\bibitem{P96}
J.~Pitman.
\newblock Random discrete distributions invariant under size-biased permutation.
\newblock {\em Adv. in Appl. Probab.}, 28(2):525--539, 1996.

\bibitem{P96b}
J.~Pitman.
\newblock Some developments of the {B}lackwell-{M}ac{Q}ueen urn scheme.
\newblock {\em Statistics, probability and game theory}, volume 30 of {\em IMS Lecture Notes Monogr. Ser.}.
\newblock Inst. Math. Statist., Hayward, CA, pp. 245--267, 1996.


\bibitem{P92}
J.~Pitman.
\newblock Poisson-{K}ingman partitions.
\newblock {\em Statistics and science: a Festschrift for Terry Speed}, volume 40 of {\em IMS Lecture Notes Monogr. Ser.}.
\newblock Inst. Math. Statist., Beachwood, OH, pp. 1--34, 2003.

\bibitem{P02}
J.~Pitman.
\newblock {\em Combinatorial stochastic processes}, volume 1875 of {\em Lecture
  Notes in Mathematics}.
\newblock Springer-Verlag, Berlin Heidelberg, 2006.
\newblock Lecture notes from \'{E}cole d'\'{e}t\'{e} de Probabilit\'{e}s de
  Saint-Flour XXXII - 2002, {\em ed.} J. Picard.


\bibitem{WM97}
S.~Walker and P.~Muliere.
\newblock Beta-{S}tacy processes and a generalization of the {P}\'olya-urn
  scheme.
\newblock {\em Ann. Statist.}, 25(4):1762--1780, 1997.

\bibitem{WM99}
S.~Walker and P.~Muliere.
\newblock A characterization of a neutral to the right prior via an extension of {J}ohnson's sufficientness postulate.
\newblock {\em Ann. Statist.}, 27(2):589--599, 1999.

\bibitem{W84}
G.~Watterson.
\newblock Lines of descent and the coalescent.
\newblock {\em Theor. Popul. Biol.}, 26(1):77--92, 1984.

\bibitem{Z}
S. Zabell
\newblock Predicting the unpredictable,
\newblock {\em Synthese}, 90(2):205--232, 1992.


\end{thebibliography}


\end{document}